\documentclass[%
 reprint,
 amsmath,amssymb,
 aps,
]{revtex4-2}

\usepackage{graphicx}
\usepackage{dcolumn}
\usepackage{bm}
\usepackage{amsmath}
\usepackage{subcaption}
\usepackage{float}
\usepackage{hyperref}
\usepackage{multirow}
\usepackage{xcolor}
\usepackage[mathlines]{lineno}% 
\usepackage{fancyhdr}
\usepackage{lipsum} 
\makeatletter
\newcommand*{\rom}[1]{\expandafter\@slowromancap\romannumeral #1@}

\begin{document}

\title{Hovering Flight in Flapping Insects and Hummingbirds: A Natural Real-Time and Stable Extremum Seeking Feedback System}

\author{Ahmed A. Elgohary}
 \altaffiliation{PhD student at the Department of Aerospace Engineering and Engineering Mechanics, University of Cincinnati, Ohio, USA}
 \email{elgohaam@mail.uc.edu}
 
\author{Sameh A. Eisa}
\altaffiliation{Assistant professor at the Department of Aerospace
Engineering and Engineering Mechanics, University of Cincinnati, Ohio,USA}
\email{eisash@ucmail.uc.edu}

\begin{abstract}
In this paper, we take an initial and novel step toward characterizing the physics of the hovering phenomenon in flapping insects and hummingbirds as a new class of extremum seeking (ES) feedback systems. By characterizing hovering flight in insects and hummingbirds as a natural hovering ES system, we achieve: (1) very simple, (2) stable, (3) model-free, and (4) real-time hovering. \textcolor{black}{More importantly, our hovering ES characterization only needs the natural oscillations of the wing as the ES input. That is, unlike other control techniques in the literature, the natural hovering ES system only needs the natural flapping action built in the system, and feedback of local sensations (measurements) related to the altitude where the insect seeks to stabilize itself. Said ES characterization, can become an important initial step in starting a new line of research that may succeed in resolving the long-standing gap between model-based control theory and the biologically observed mechanisms that stabilize hovering flight.}
% In this paper, we \textcolor{black}{take an initial and novel step towards characterizing} the physical phenomenon of hovering flight as an extremum seeking (ES) feedback system. \textcolor{black}{We anticipate that} said novel characterization \textcolor{black}{may start a new line of research that can potentially} solves all the puzzle pieces of hovering flight that existed for decades in previous literature. Is hovering flight stable? If so, what is the control mechanism utilized by insects/hummingbirds to achieve stable hovering? If such a mechanism exists, does it fit the biological constraints that insects/hummingbirds have limited computational abilities? Does it fit the experimental biology narrative that insects/hummingbirds rely mainly on their sensation to stabilize hovering? Our ES characterization and analysis provide for the first time a simple, model-free, real-time, stable feedback system of hovering. Consistent with natural observations and biological experiments, hovering via ES is simply achievable by the natural oscillations of the wing angle and measuring (sensing) altitude or \textcolor{black}{acceleration}. 
We provide simulation trials, \textcolor{black}{including comparisons with some approaches from literature,} to demonstrate the effectiveness and robustness of our results. \textcolor{black}{We used literature data for} hawkmoth, cranefly, bumblebee, dragonfly, hoverfly, and a hummingbird.

\end{abstract}

%\tableofcontents
\maketitle
\section{Introduction}
Flapping flight is a fascinating physical phenomenon that captivated the interests of different communities such as physicists, biologists, and aerospace/control engineers. For decades, flapping flight has been a rich research area with a highly complex nature, involving multi-body, nonlinear, and non-autonomous dynamics \cite{taha_review,phan2019insect}. Moreover, the mechanics of flapping-wing flight impose a unique challenge on researchers in many fields, as conventional aerodynamic principles developed for fixed-wing aircraft are quite limited in characterizing flapping-wing systems \cite{taha_review,xuan2020recent}. This is especially true when it comes to the flapping flight of insects and hummingbirds due to the ability of said organisms to \textit{hover}. Hovering flight is a unique, puzzling physical phenomenon for biologists and engineers. It was once deemed impossible under traditional aerodynamic theory, as the required hovering lift coefficients for balance are significantly higher than what is achievable through conventional aerodynamics \cite{weis1972energetics,norberg1975hovering,dudley1990mechanics,ellington1995unsteady,145,taha_review,sun2014insect}. Over time, unconventional lift mechanisms that can physically enable hovering flight has been discovered, particularly the concept of the leading-edge vortex, introduced in \cite{nature1996LEV} and later expanded upon by the authors in \cite{dickinson1999wing}, which provides the necessary lift for insect hovering flight. Consequently, significant progress has been made in both theoretical and experimental research to better understand, analyze, and mimic hovering flight \cite{taha_review,phan2019insect,xuan2020recent}. On the theoretical side, researchers tried to derive and advance the modeling aspects of hovering flight (on the body and aerodynamic levels) to analyze properties such as stability. They also aimed to introduce control mechanisms/techniques that can explain stable hovering and enable the application of hovering in robotics. On the experimental side, most efforts are aimed at studying flapping insects and the relation between their sensation and the ability to hover.
% \section{Introduction}

% Flapping flight is a fascinating physical phenomenon that captivated the interests of physicists, biologists, and aerospace/control engineers. For decades, flapping flight has been a rich research area with a highly complex nature, involving multi-body, nonlinear, and non-autonomous dynamics \cite{taha_review,phan2019insect}. Hovering flight is a unique, puzzling physical phenomenon for biologists and engineers: it was once deemed impossible under traditional aerodynamic theory, as the required hovering lift coefficients for balance are significantly higher than what is achievable through conventional aerodynamics \cite{weis1972energetics,norberg1975hovering,dudley1990mechanics,ellington1995unsteady,145,taha_review,sun2014insect}. Over time, unconventional lift mechanisms that can physically enable hovering flight have been discovered, particularly the concept of the leading-edge vortex, introduced in \cite{nature1996LEV} and later expanded upon in \cite{dickinson1999wing}, which provides the necessary lift for insect hovering flight. **Consequently, significant progress has been made in both theoretical and experimental research to better understand, analyze, and mimic hovering flight** \cite{taha_review,phan2019insect,xuan2020recent}.

\smallskip
The long-held view and consensus in the literature is that hovering flapping flight is inherently unstable \textcolor{black}{in the open-loop sense} and active control may be required for stabilization \cite{taha_review,sun2014insect,taylor2003dynamic,taylor2005nonlinear,zhang2012lateral,xu2013lateral,DelftDissertation2014,feedback25ristroph2010discovering}. A recent study \cite{lyu2022dynamic} extended stability analysis to a broader range of hovering insects has reached the same conclusion as well. Recently, the authors in \cite{taha2015need,taha2020vibrational} argued that concluding instability is due to a deficiency in first-order averaging and that hovering flight may be shown stable—including the pitching motion in an open-loop sense if higher-order averaging is used; however, they acknowledged that feedback control may be required for other insects to stabilize their hovering flight. Natural and biological observations/experiments are in consensus that flapping insects rely heavily on their sensation, which means it is highly plausible that they use some feedback control mechanism to stabilize their hovering in real time \cite{fuller2014flying,sensorroyalristroph2013active,taylor2003dynamic,feedback25ristroph2010discovering,feedback26cheng2011mechanics,sensors1taylor2007sensory}. 
% \textcolor{black}{The authors in \cite{feedback25ristroph2010discovering,sensorroyalristroph2013active} examined how freely flying fruit flies react in real time to external disturbances affecting their yaw and pitch: small ferromagnetic materials were attached to the flies, allowing for reorientation via a magnetic field, and observations indicated that the flies adapted their wing kinematics in real-time.} \textcolor{black}{Thus, this can imply that, in the absence of passive damping, active sensorimotor
% feedback is important for stabilization.}

\smallskip
\textcolor{black}{From a control-theoretic perspective, stabilized hovering is achieved through sensory feedback \cite{taha2020vibrational,IEEETransaction,taylor2003dynamic,feedback25ristroph2010discovering,rapp2020spiking,fuller2014flying}. This understanding has motivated many researchers to develop feedback-based control techniques aimed at enabling stable hovering. Consequently, a variety of model-dependent control strategies have been proposed in the literature \cite{taha_review,sun2014insect}}. Standard approaches such as Linear Quadratic Regulator (LQR), Proportional-Integral-Derivative (PID) control, and other nonlinear designs rely on dynamic and aerodynamic models that require knowledge of the equilibrium hovering point \emph{a priori} \cite{24,Delft_Conference,serrani2010robustAIAA,doman2009dynamics}. In several studies, aerodynamic models necessary for control were based on computational fluid dynamics simulations, which are inherently non-real-time \cite{oppenheimer2009dynamics,Delft_Conference,sun2004computational,hedrick2015recent}. \textcolor{black}{It is worth mentioning that, although PID and its variants (PI and PD) are real-time controllers, even the PI scheme that reproduces yaw recovery in fruit flies \cite{feedback25ristroph2010discovering} — as the authors themselves state — requires prior estimation of morphological and aerodynamic parameters and subsequent hand-tuning of gain parameters, with all parameters varying by $\pm 15\%$.
 Moreover, none of the control methods, techniques, or mechanisms discussed above leverage the inherent high-frequency wing-flapping perturbation action, which can serve as a built-in probing signal.}

% for extremum-seeking. Very simple, real-time, model-free control for hovering in flapping insects cannot be found in literature prior to this paper, to the best of our knowledge, and this gap motivates the contributions of the present work.

% \textcolor{red}{As noted in literature, model-based control designs have their downsides since they inherently depend on the accuracy of the dynamic/aerodynamic models. Increasing their accuracy is quite challenging computationally; we leave the reader with this quotation from \cite{IEEETransaction} which expands on this point:``Although, at present, some numerical simulations of unsteady insect flight aerodynamics based on the finite element solution of the Navier–Stokes equations give accurate results for the estimated aerodynamics forces \cite{IEEE27,IEEE28}, their implementation is unsuitable for control purposes since they require several hours of processing for simulating a single wingbeat, even on multiprocessor computers."}

\smallskip
\noindent\textbf{Motivation.}
The consensus that hovering flight is open-loop unstable and needs feedback control for stabilization has motivated researchers to introduce a variety of feedback control techniques to achieve stable hovering; however, available control methods can be model-dependent, require knowledge of the equilibrium hovering point \emph{a priori}, and can potentially be computationally expensive to implement in real-time due to tuning challenges and/or high need of computational power. \textcolor{black}{Moreover, none of these controllers leverage the inherent high-frequency wing-flapping perturbation action built into the system. Hence, we are motivated to propose a control method that can frame the hovering problem while using the built-in high-frequency wing perturbation action in a natural way as the probing signal for the controller itself. For that, we turn our attention to Extremum Seeking (ES) control systems \cite{ariyur2003real,scheinker2017model,scheinker2024100}, which utilize perturbation actions as the probing signal for the controller to derive a given dynamical system towards the extremum of unknown objective functions in a model-free, real-time, sensation-based fashion. Moreover, ES feedback control systems have achieved significant successes in characterizing and mimicking challenging behaviors in biological organisms immersed in fluid structures — see, for example, fish \cite{ESCfishcochran2009source} and microswimmers \cite{krstic2008extremum,abdelgalil2022sea}. Recently, our group introduced ES characterization in flying organisms and succeeded in decoding the energy-efficient soaring strategy of birds \cite{pokhrel2022novel,eisa2023analyzing}.}  

\noindent\textbf{Contribution.} In this paper, we take an initial and novel step toward characterizing the physics of the hovering phenomenon in flapping insects and hummingbirds as a new class of ES feedback systems. By characterizing hovering flight in insects and hummingbirds as a natural hovering ES system, we \textcolor{black}{achieve}: (1) very simple, (2) stable, (3) model-free, and (4) real-time hovering. More importantly, our hovering ES characterization only needs the natural oscillations of the wing (i.e., the natural flapping action, built in the system, unlike other control techniques in the literature) as the ES input, and feedback of \textcolor{black}{local} sensations (measurements) related to the altitude where the insect seeks to stabilize itself. Said ES characterization, can become an important initial step in starting a new line of research that may succeed in resolving the long-standing gap between model-based control theory and the biologically observed mechanisms that stabilize hovering flight. It also may provide model-free, real-time, computationally cheap, bio-mimicry possibilities for flapping robots.

In this paper, we test our proposed hovering ES (see Section \ref{Sec: Natural hovering flight as extremum seeking system}) \textcolor{black}{using literature data} for five flapping insects (hawkmoth, cranefly, bumblebee, dragonfly, hoverfly) and a hummingbird using simplified models (see Section \ref{sec:Modelling}) representing the body dynamics and aerodynamics of flapping insects/hummingbirds. Said testing is for the proof of concept as ES systems work in a model-free fashion \cite{ariyur2003real,scheinker2017model,scheinker2024100}. \textcolor{black}{We also compare our ES approach with open-loop setting and a PID control approache, \textcolor{black}{including optimal-gains PID controllers}.} Testing and simulation results, \textcolor{black}{including with feedback noise and delay}, indicate success in stabilizing hovering for the aforementioned insects and hummingbird \textcolor{black}{based on literature data}. Following that (in Section \ref{sec:Stability}) we also provide stability analysis using sophisticated tools from differential geometric control theory \cite{bullo2002averaging,pokhrel2023higher,elgohary2025extremum}. 
% for the insects/humminbird we tested in Section \ref{Sec: Natural hovering flight as extremum seeking system}. 
Stability results indicate stable hovering for all tested insects and hummingbird.
% ; stability has never been concluded in the literature for cranefly, bumblebee, dragonfly, hoverfly, and hummingbirds. 
Finally, we conclude the paper in Section \ref{sec:conclusion}.

\section{Modelling of Flapping Insects and Hummingbirds: Body Dynamics and Aerodynamics}
\label{sec:Modelling}

Even though ES is a model-free system as mentioned earlier, we use the model provided in this Section for testing, simulation, and verification (i.e., as a proof of concept). Now, we lay the groundwork for the model. 
% In this section, we briefly provide the aerodynamics and dynamic flapping models used for our simulation results.  
% Remark: recall from the contribution section that our controller is model-free; however, in this study, we employed a flapping model for the purpose of analysis and to validate the simulation results.

Building models for flapping flight dynamics is particularly challenging due to the oscillatory motion of the wings relative to the body \cite{taha_review}. Hence, models of flapping flight are usually represented as multi-body, nonlinear time-periodic (NLTP) systems \cite{taha_review,taha2014longitudinalmodelguidance,taha2015need}. Moreover, these models operate as multi-scale dynamical systems, given that the body flight dynamics and wing flapping dynamics operate on distinct time scales \cite{hassan2018combinedguidancejournal,taha2020vibrational}; this separation in time scales between the body (slow) and wings (fast) can be visually observed in nature. 
% Moreover, these models operate as multi-scale dynamical systems, given that the body flight dynamics and wing flapping dynamics operate on distinct time scales \cite{hassan2018combinedguidancejournal,taha2020vibrational}.

\textcolor{black}{There are} two primary approximations that have been considered in flapping flight modeling, namely physical approximations and mathematical approximations. We only revisit the physical approximations since they can affect the accuracy of the model itself (note that mathematical approximations are merely for stability analysis, so they do not affect the modeling accuracy when simulating the flapping system). Recently, there have been some developments \cite{taha2015need,taha2013unsteady,141wang2004unsteady,1andersen2005unsteady,taha2016geometriclong,taha2014longitudinalmodelguidance,hassan2018combinedguidancejournal,hassan2017combinedscitech} that we think are worth consideration in our modeling.
% As previously discussed in section IA, two primary approximations have been adopted for stability analysis and control studies of flapping flight. However, those approximations become refuted and should be reconsidered for a more accurate analysis of flapping flight dynamics. 
In \cite{taha2016geometriclong}, the authors demonstrated that the interaction between wing flapping dynamics and body flight dynamics can lead to a negative lift mechanism, significantly affecting overall flight performance. Additionally, on the aerodynamic side, unsteady aerodynamics models \cite{taha2013unsteady,141wang2004unsteady,1andersen2005unsteady} for flapping have gained traction as they include factors like leading-edge vortex, rotational dynamics, viscous and added mass effects, and wake capture, all of which contribute to a more accurate representation of aerodynamic loads.
% numerous studies have shown that the aerodynamic forces in flapping flight are primarily driven by the leading-edge vortex and rotational effects \cite{Delft_journal,de2024bio,lentink2009rotational,sane2002aerodynamic}. As a result,   
% Recently, interest in aerodynamic models of flapping flight has grown to better understand the complex flow fields involved. In particular, the unsteady nature of aerodynamics models, where factors like leading-edge vortex, rotational dynamics, viscous and added mass effects, and wake capture contribute to the overall aerodynamic loads \cite{taha2013unsteady,141wang2004unsteady,1andersen2005unsteady}. Developing models that can capture these complex, nonlinear, and unsteady phenomena without excessive computational cost remains a major challenge. Numerous studies have shown that the aerodynamic forces in flapping flight are primarily driven by the leading-edge vortex and rotational effects, as demonstrated in \cite{Delft_journal,de2024bio,lentink2009rotational,sane2002aerodynamic}
\begin{figure}[h]
    \centering
\includegraphics[width=1\linewidth]{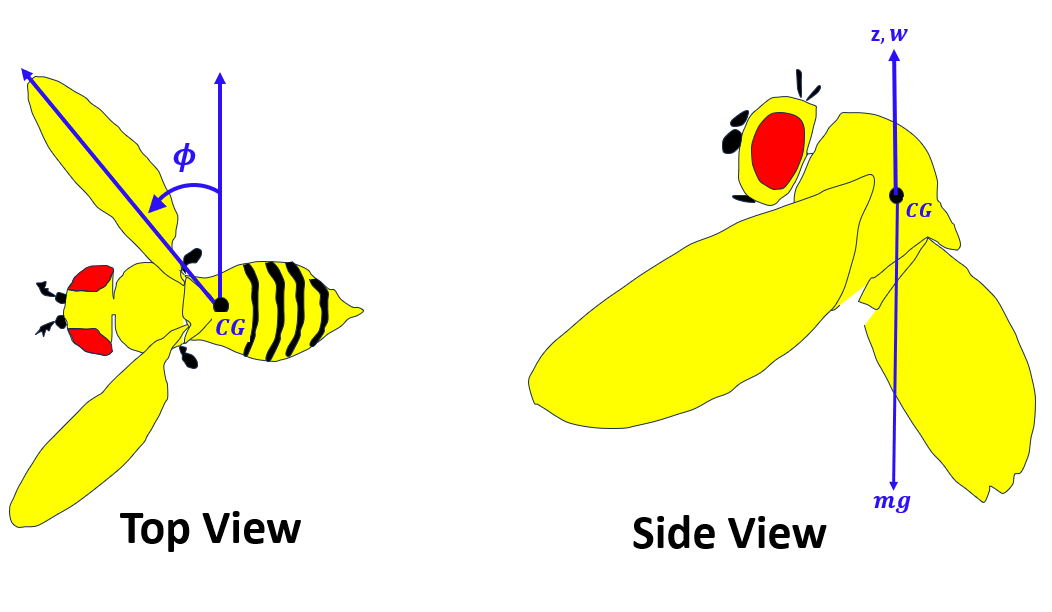}
    \caption{Schematic body diagram of the flapping insect system which moves vertically via the vertical velocity $w$ (to gain or lose altitude $z$) and the flapping angle $\phi$.}
    \label{Body_digram}
\end{figure}

\begin{figure}[h]
    \centering
\includegraphics[width=0.9\linewidth]{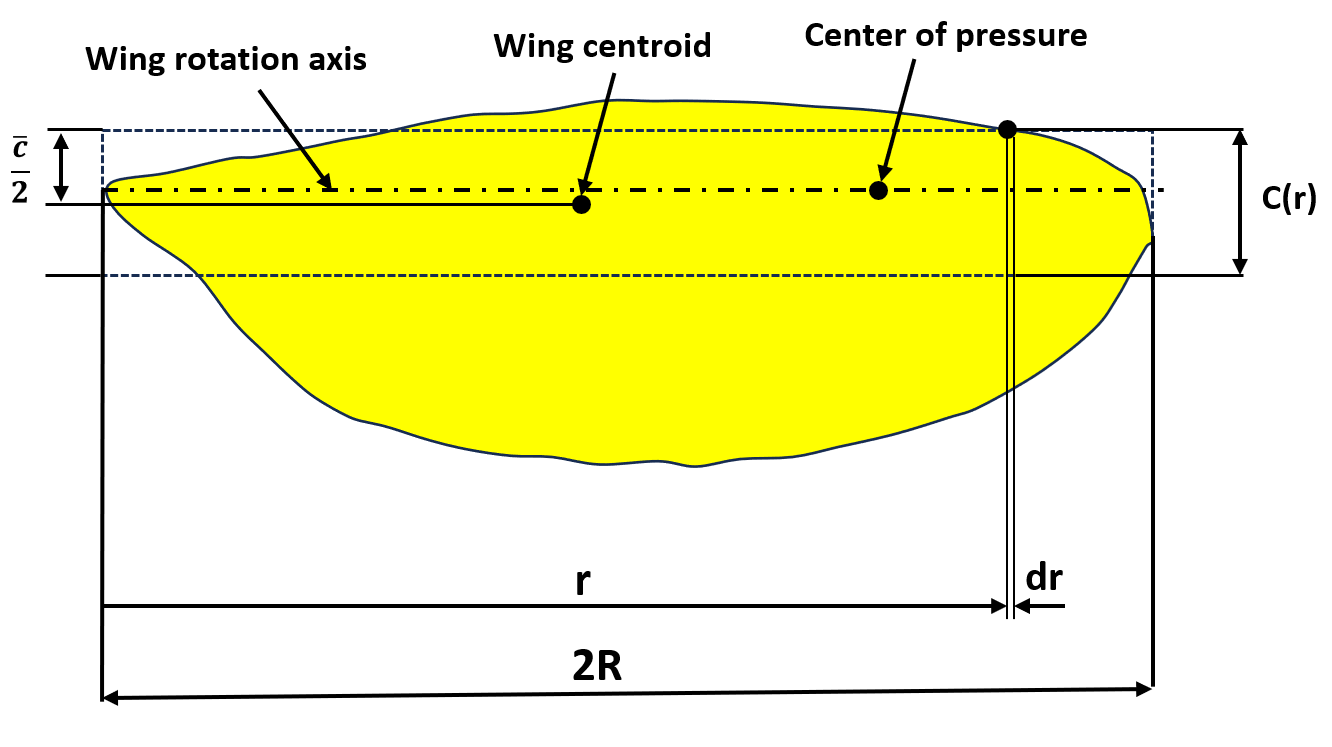}
    \caption{Sketch of main wing parameters.}
    \label{fig:wing_parameters}
\end{figure}

Given all the above, the flapping flight model in this paper will be treated as a multi-body problem, accounting for both the wing and body as rigid bodies. Furthermore, key aerodynamic features, such as the leading-edge vortex and rotational effects, should be incorporated into the aerodynamic models. As a result, we turn our attention to \cite{taha2016geometriclong,taha2014longitudinalmodelguidance} where a flapping flight dynamics model near hover was analyzed using multi-rigid-body equations, treating both the wing and body as rigid bodies. In addition, the authors in \cite{taha2016geometriclong,taha2014longitudinalmodelguidance} also considered a simplified analytical aerodynamic model which is derived to account for key aerodynamic effects such as the leading-edge vortex and rotational effects.
\begin{table*}
\centering
 \caption{The morphological parameters for the five insects and hummingbird we study in this paper.}
\label{tab:morphological_parameters}
\resizebox{\textwidth}{!}{%
\begin{tabular}{lccccccccc}
\hline
\textbf{Insect} & \textbf{f (Hz)} & \boldmath$\Phi^\circ$ & \textbf{S (mm\(^2\))} & \textbf{R (mm)} & \boldmath$\bar{c}$ \textbf{(mm)} & \boldmath$\hat{r}_1$ & \boldmath$\hat{r}_2$ & \textbf{m (mg)} \\
\hline
Hawkmoth    & 26.3  & 60.5 & 947.8 & 51.9 & 18.3 & 0.440 & 0.525 & 1648  \\
Cranefly    & 45.5  & 61.5 & 30.2  & 12.7 & 2.38 & 0.554 & 0.601 & 11.4  \\
Bumblebee   & 155   & 58.0 & 54.9  & 13.2 & 4.02 & 0.490 & 0.550 & 175   \\
Dragonfly   & 157   & 54.5 & 36.9  & 11.4 & 3.19 & 0.481 & 0.543 & 68.4  \\
Hoverfly    & 160   & 45.0 & 20.5  & 9.3  & 2.20 & 0.516 & 0.570 & 27.3  \\
Hummingbird & 48    & 70   & 611   & 48   & 12.7 & 0.428 & 0.492 & 4320  \\
\hline
\end{tabular}%
}
\end{table*}
% In \cite{taha2016geometriclong,taha2014longitudinalmodelguidance}, a flight dynamics model for flapping flight near hover was analyzed using multi-rigid-body equations, treating both the wing and body as rigid bodies. In addition, a simplified analytical aerodynamic model was derived to account for key aerodynamic effects such as the leading-edge vortex and rotational effects. 
To expand further, in \cite{taha2016geometriclong} longitudinal equations of motion for the body-wing system were derived using the principle of virtual power. Building upon the longitudinal model in \cite{taha2016geometriclong}, the same authors  \cite{hassan2018combinedguidancejournal,hassan2017combinedscitech} developed a simple, yet accurate two-degree-of-freedom (2-DOF) dynamics model that captures the essential dynamics of the system, characterized by vertical velocity $w$ and the wing flapping angle $\phi$. Hence, the body moves vertically (i.e., studying hovering based on altitude).
% developed a two-degree-of-freedom (2-DOF) dynamics model, in which the body is constrained to move along vertical rails, characterized by vertical velocity $w$ and the wing flapping angle $\phi$. 
In this paper, we utilize this 2-DOF model by \cite{hassan2018combinedguidancejournal,hassan2017combinedscitech} for analyzing the flapping flight dynamics of hovering. 
% as it provides a simplified but effective framework to capture the essential dynamics of the system. 
The governing equations for the 2-DOF model are as follows (see Figure \ref{Body_digram}): 
\begin{equation}
    \begin{aligned}
        \dot{w} &= g - k_{d1} |\dot{\phi}| w - k_L \dot{\phi}^2, \\
        \ddot{\phi} &= -k_{d2} |\dot{\phi}| \dot{\phi} - k_{d3} w \dot{\phi} + \frac{1}{I_F} \tau,
    \end{aligned}
    \label{model_equation}
\end{equation}
where \( g \) is the gravitational acceleration, \( I_F \) is the wing flapping moment of inertia given in \cite{taha2016geometriclong}, and \(\tau \) is the input control torque. The coefficients \( k_{d1}, k_L, k_{d2}, k_{d3} \) are defined as follows:
\begin{equation}
\begin{aligned}
    k_{d1} &= \frac{\rho C_{L\alpha} I_{11} \cos^2 \alpha_m}{2m}, \\
    k_L &= \frac{\rho C_{L\alpha} I_{21} \sin \alpha_m \cos \alpha_m}{2m}, \\
    k_{d2} &= \frac{\rho C_{L\alpha} I_{31} \sin^2 \alpha_m}{I_F}, \\
    k_{d3} &= \frac{\rho C_{L\alpha} I_{21} \sin \alpha_m \cos \alpha_m}{I_F},
\end{aligned}
\label{eq:coefficients}
\end{equation}
where $\rho$ represents the air density, $\alpha_m$ is the mean angle of attack sustained throughout the flapping hovering. The total mass of the insect, including both wing and body mass, is represented by $m $, and the constants $I_{mn}$ are calculated based on the chord distribution $c(r)$ along the wing, given in \cite{taha2016geometriclong} by the integral $I_{mn} = 2 \int_0^R r^m c^n(r) dr$. Note that $C_{L_\alpha}$ is the lift curve slope of the wing as defined in \cite{taha2014longitudinalmodelguidance}:
\begin{equation}
C_{L_\alpha} = \frac{\pi AR}{  1 + \sqrt{\left(\frac{\pi AR}{a_0}\right)^2 + 1}},
\label{eq:cl_alpha}
\end{equation}
where $AR$ is the wing aspect ratio given as $\frac{R^2}{S}$ and $a_0$ is the lift curve slope of the two-dimensional (2-D) airfoil. The chord distribution $c(r)$ along the wing is defined in \cite{taha2015need} as:
\begin{equation}
    c(r) = \frac{\bar{c}}{\beta} \left( \frac{r}{R} \right)^{\kappa - 1} \left( 1 - \frac{r}{R} \right)^{\gamma - 1},
\end{equation}
where
\[
\kappa = \hat{r}_1 \left[ \frac{\hat{r}_1 (1 - \hat{r}_1)}{\hat{r}_2^2 - \hat{r}_1^2} - 1 \right], \quad 
\gamma = (1 - \hat{r}_1) \left[ \frac{\hat{r}_1 (1 - \hat{r}_1)}{\hat{r}_2^2 - \hat{r}_1^2} - 1 \right],
\]
and
\[
\beta = \int_0^1 \hat{r}^{\kappa - 1} (1 - \hat{r})^{\gamma - 1} \, d\hat{r},
\]
where $\hat{r}_1$ and $\hat{r}_2$ referees to the moments of the wing chord distribution defined in \cite{taha2015need} as:
\begin{equation}
    \hat{r}_k^k = \frac{\int_0^R r^k c(r) \, dr}{S R^k}.
\end{equation}
$I_F$ is the moment of inertia of the flapping wing given in \cite{taha2016geometriclong} as:
\begin{equation}
I_F = I_x \sin^2 \alpha_m + I_z \cos^2 \alpha_m,
\label{eq:inertia_formula}
\end{equation}
where the wing inertial properties  $I_x$ and $I_z$ are estimated as follows:
\begin{equation}
\begin{aligned}
    I_x &= 2 \int_0^R m' r^2 c(r) dr = m' I_{21}, \\
    I_y &= 2 \int_0^R m' \hat{d}^2 c^3(r) dr = m' \hat{d}^2 I_{03}, \\
    I_z &= I_x + I_y,
\end{aligned}
\end{equation}
where $m'$ represents the distribution of the wing's areal mass, defined as $\frac{m_w}{2S}$, where $m_w$ is the wing mass \cite{wu2009hovering}. Additionally, $\hat{d}$ is the distance from the wing's hinge line to its center of gravity (CG), normalized by the chord length \cite{ellington1984aerodynamics}.
The main wing parameters are shown in Figure \ref{fig:wing_parameters}. 

For the convenience of the reader, we collected all needed parameters (and made all needed calculations for some parameters) for implementation and simulations, and provide them here in this Section. The morphological parameter values for the five insects we study in next Section and the hummingbird \cite{taha2015need,Delft_Conference} are provided in Table \ref{tab:morphological_parameters}, where $f$ is the flapping frequency, $\Phi^\circ$ is the amplitude of the flapping angle, $S$ is the surface area of one wing, $R$ is the wing radius, $\bar{c}$ is the mean chord length, $m$ is the body mass, and $I_y$ is the body moment of inertia around the y-axis of the body. The mean angle of attack, $\alpha_m$, and the mass of a single wing were given in \cite{ellington1984aerodynamics} for the hoverfly, dragonfly, and cranefly insects, in \cite{145} for the hawkmoth, in \cite{dudley1990mechanicsbumblebees} for the bumblebee, and in \cite{Delft_Conference} for the hummingbird. In Table \ref{tab:Model parameters}, the evaluated values of the parameters \(k_{d1}\), \(k_L\), \(k_{d2}\), \(k_{d3}\), and \(I_F\) for each insect and the hummingbird are provided.
\begin{table}[h]
\centering
\caption{Calculated model parameters for the five insects and hummingbird we study in this paper.}
\label{tab:Model parameters}
\begin{tabular}{lccccc}
\hline
\textbf{Insect} & \boldmath$k_{d1}$ & \boldmath$k_L$ & \boldmath$k_{d2}$ & \boldmath$k_{d3}$ & \boldmath$I_F$ \\
\hline
Hawkmoth    & 0.0354 & 6.216e-04 & 0.3492 & 17.3331 & 1.3179-07 \\
Cranefly    & 0.0787 & 2.647e-04 & 0.4237 & 114.9124 & 5.703e-11 \\
Bumblebee   & 0.0072 & 2.204e-05 & 0.2826 & 82.5021 & 9.453e-11 \\
Dragonfly   & 0.0115 & 2.548e-05 & 0.1127 & 45.1702 & 7.956e-11 \\
Hoverfly    & 0.0143 & 3.715e-05 & 0.3099 & 106.9442 & 1.945e-11 \\
Hummingbird & 0.099 & 1.549e-04 & 0.5515 & 30.6250 & 4.452e-08 \\
\hline
\end{tabular}
\end{table}

Finally, the dynamics of the system \ref{model_equation} can be reformulated into state-space representation to include also the altitude $z$ and flapping angle $\phi$ as:
\begin{equation}
\label{eq:Model_state_space}
\frac{d}{dt}
\begin{bmatrix}
z \\
\phi \\
w \\
\dot{\phi}
\end{bmatrix}
=
\begin{bmatrix}
w \\
\dot{\phi} \\
g - k_{d1} |\dot{\phi}| w - k_L \dot{\phi}^2 \\
-k_{d2} |\dot{\phi}| \dot{\phi} - k_{d3} w \dot{\phi}
\end{bmatrix}
+ 
\begin{bmatrix}
0 \\
0 \\
0 \\
\frac{1}{I_F}
\end{bmatrix}
\tau.
\end{equation}

\textcolor{black}{\textbf{Remark.}
It is to be noted that the model in (\ref{eq:Model_state_space}) is a reduced and simplified model for hovering. Nevertheless, as long as this reduced model is relatively accurate in capturing the essential qualitative features of hovering, it is sufficient for testing the ES framework. This is due to the fact that ES systems are inherently model–free and tolerant of large model uncertainties and/or variations \cite{scheinker2024100}; in our case, this will be further illustrated in Section III.D. That is, successful performance and stabilization in hovering for the simplified model adopted in this work, suggest that the same ES framework will remain effective under the richer, more variable dynamics of real insects. In other words, the ES framework is robust to modeling errors.}

\section{Natural Hovering flight as an extremum seeking system}
\label{Sec: Natural hovering flight as extremum seeking system}
Extremum seeking (ES) systems \cite{ariyur2003real,scheinker2017model,scheinker2024100} are themselves control systems (called sometimes ESC). They steer/drive a dynamical system to the extremum point (maximum/minimum) of an objective function that can be unknown, expression-wise. ES systems have been powerful and attractive for many reasons, including that they are (1) very simple structurally and in implementation, (2) stable, (3) model-free, and (4) real-time. In fact ES systems have found many applications in real-time dynamic optimization and control systems in many fields (for more details, the reader can refer to \cite{scheinker2024100}). With relevance to the topic of this paper, ES systems have had significant successes in characterizing and mimicking some optimized behaviors in biological organisms such as fish \cite{ESCfishcochran2009source}, micorswimmers \cite{krstic2008extremum,abdelgalil2022sea}, and soaring birds \cite{pokhrel2022novel,eisa2023analyzing}.  

ES systems \cite{krstic2000stability,ariyur2003real,scheinker2024100}, generally speaking, can be considered as input-output systems. The input is a perturbation (modulation) action/signal, which causes variation in the system dynamics. This translates into variations (corresponding to input perturbations/modulation) in the objective function, which is considered the output in ES systems and the only feedback required. Recently, inspired by the flapping mechanics of insects, a class of ES systems was proposed \cite{elgohary2025extremum}, which is relevant to systems possessing mechanical structure (e.g., systems described by Newton's second law or Euler-Lagrange formulation) with generalized forces that are quadratic in velocities; this matches, and can be applied to, the system representing hovering at hand in (\ref{eq:Model_state_space}).    
% As long as the measurements of the objective function are accessible, the ES system in fact does not need even the system dynamics; hence, the model-free nature of ES systems.
\textcolor{black}{Unlike most control techniques, in ES systems, even the objective function, including its maximum/minimum point, can be unknown. All that the ES system needs is access to its \textcolor{black}{local} measurements. Unfamiliar readers with why this model-free nature of ES (including dealing with unknown objective functions) is particularly useful for the paper at hand, and is a near-unique property among control techniques, may refer to these two papers regarding model-free, real-time light-source seeking \cite{ECC2024,unicycle_bio_inspiration} by ES systems -- see also this video \cite{Unicycle_multiple_positions}. In these papers (and the video), the light source, which represents the maximum of the objective function, that is, the light signal itself, is unknown position-wise; the ES system steers the robot to the light source without any information about the position of the source, but by simply measuring (sensing) light.}   
% ES systems andThe reader can refer to model-free light source seeking by ES systems here \cite{ECC2024,unicycle_bio_inspiration} These measurements are the only feedback the ES system needs to update its input, and start the cycle anew. 
In this paper, we characterize and treat hovering flight as an ES system (i.e., ES is the model and control simultaneously for hovering). First, we expand on the philosophy/reasoning behind our hypothesis before moving into our attempt to provide a proof of concept.  

Hovering insects/hummingbirds \cite{taha_review,phan2019insect,xuan2020recent} utilize flapping as their action. To be more precise, and to expand the word ``flapping" in physical terms, they conduct very high oscillatory motion of their flapping wing angle $\phi$. In kinematic terms, this is usually translated into a high amplitude sinusoidal function. For the flapping model we use in this paper (see \eqref{eq:Model_state_space}), said oscillatory motion is associated with the torque $\tau$, which is responsible for changing the flapping angle rate $\dot{\phi}$. Thus, perturbations in $\tau$ is the modulation step. Due to perturbation/modulation in $\tau$, $\dot{\phi}$ is updated, causing variation in $\phi$. As a consequence, the dynamic system (including the insect's body) will be varied as well, resulting in variation in any measurements (sensation by the insect) that depend on the state variables of the system dynamics. \textcolor{black}{Let us clarify this further using the light source seeking analogy we discussed earlier (including citations to these papers \cite{ECC2024,unicycle_bio_inspiration} and this video \cite{Unicycle_multiple_positions}). Let us assume that a flapping insect (e.g., a moth) is stabilizing its altitude (i.e., hovering) based on the sensation of light emitting from a source with unknown position. This means that the altitude about which the insect will stabilize corresponds to the unknown maximum light intensity. If we successfully characterize hovering flight as an ES system, this implies that via \textcolor{black}{local} measuring (sensing) altitude or altitude-related signals, the insect will be able to stabilize about the desired altitude, which corresponds to the maximum light intensity. Moreover, in source seeking applications, such as light source seeking, a perturbation signal is typically introduced into the control input as an integral part of the ES control system \cite{ECC2024,unicycle_bio_inspiration}. This perturbation enables the system to sense the objective function and adjust accordingly. However, in the natural hovering ES system, this additional perturbation is \textit{unnecessary}, as it is inherently present due to the wing flapping motion itself. The oscillatory nature of flapping wings naturally induces this perturbation, allowing the system to sense the objective function in real-time without any external modifications, even when the objective function is unknown.}

\textcolor{black}{It is worth to mention and emphasize that the light–seeking illustration is only one example of how local sensation or measurements might be enough to seek an extremum. But this local sensation (measurements) could be heat, light, chemical concentration, among many other possibilities. For example, the authors in \cite{abdelgalil2022sea} apply ES to sperm that cannot “see” the egg and rely only on chemotaxis.}
% In this sense, the mechanism can be global: mosquitoes \cite{seeking_heat} and certain beetles are drawn to warm temperatures, influencing their breeding and foraging behavior.}

\begin{figure*}[ht]
    \centering
\includegraphics[width=0.85\linewidth]{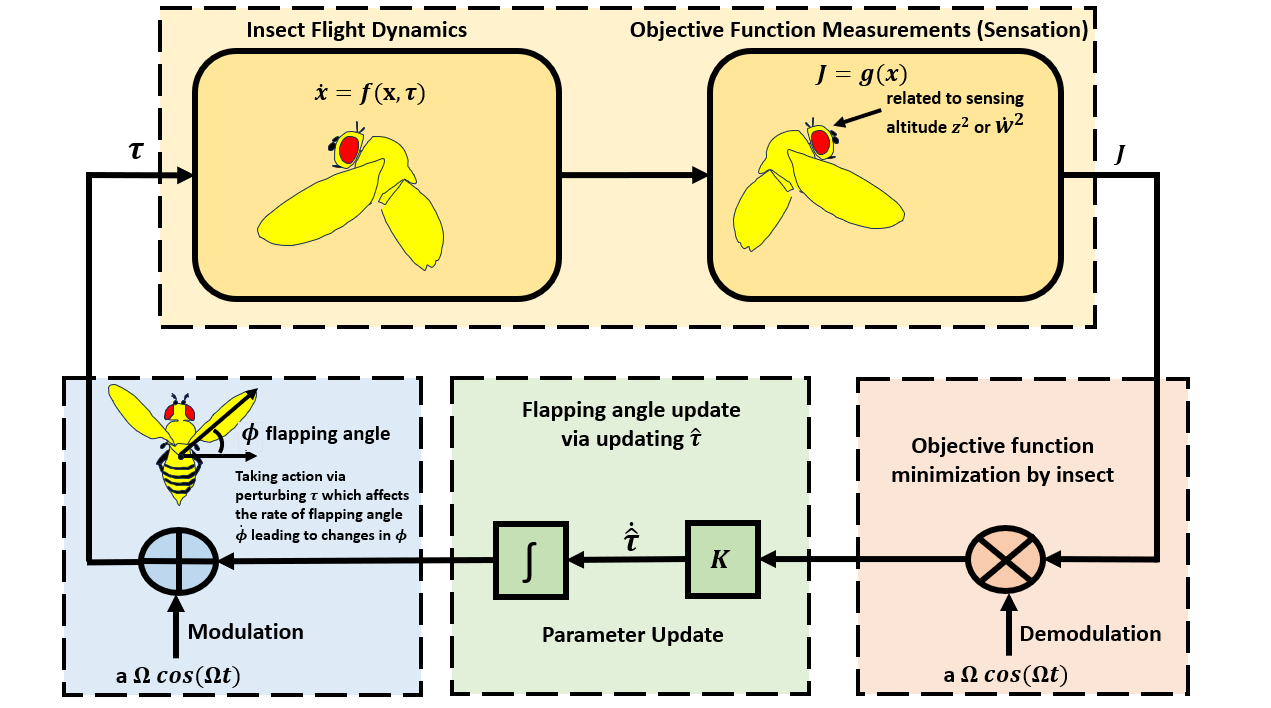}
\caption{Natural hovering extremum seeking system for flapping insects and hummingbirds. The perturbation of the torque $\tau$ is consistent of the average/mean torque estimate $\hat{\tau}$ added to it the modulation signal $a \Omega cos(\Omega t)$, i.e., $\tau = \hat{\tau} + a \Omega \cos(\Omega t)$. This perturbation in $\tau$ affects the flapping angle rate $\dot{\phi}$, causing variation in the flapping angle $\phi$. Said variations in $\phi$ is reflected in the dynamic system $\bm{\dot{x}} = f(\bm{x}, \tau)$ (the insect/hummingbird). As a result, the insect/hummingbird sensation (i.e., objective function measurements) is also varied in a domino-effect fashion. The sensation (measurements) of the objective function reflects the altitude (or altitude-related signals) or acceleration (or acceleration-related signals). Measurements of the objective function is the feedback which is processed by demodulation, providing the update for the torque average/mean $\hat{\tau}$ which takes the insect/hummingbird closer to the minimum of $J$, seeking hovering in real-time.}  \label{fig:ESC_customized_structure}
\end{figure*}

Based on the biological observations \textcolor{black}{and the light source seeking discussion above}, we rationally assume that the insect/hummingbird is able to sense their altitude (or altitude-related signals). Hence, without loss of generality, we assume the objective function $J$ to reflect the altitude difference between the current state and the desired altitude 
% or \textcolor{black}{acceleration} 
to stabilize about certain altitude. The hovering ES system aims to minimize $J$, so $J$ measurements (sensation) is demodulated, and is used to update the torque average (mean estimate) $\hat{\tau}$. Then, another cycle starts by perturbation again. Our natural hovering ES system is provided in Figure \ref{fig:ESC_customized_structure}.     

 One major difference between our proposed natural hovering ES system in Figure \ref{fig:ESC_customized_structure} and the traditional classic ES structure in \cite{krstic2000stability,ariyur2003real} is the nature of the perturbation/modulation (also the demodulation) signals. In traditional classic ES systems, the amplitude of the modulation and demodulation signals is considered small (e.g., $a\: cos(\Omega t)$ for relatively small $a$). However, in the proposed natural hovering ES system, we use high amplitude, high frequency perturbation signals in the form $a\: \Omega cos(\Omega t$) following the natural/biological observations in literature, which are reflected in the flapping model we use in \eqref{eq:Model_state_space}. As a result, our proposed natural ES system is a member of the ES classes in \cite{elgohary2025extremum}. Furthermore, the perturbation input, flapping angle torque $\tau$, is given by $\tau = \hat{\tau} + a \Omega \cos(\Omega t)$. The objective function $J$ which represents the output $J = g(x)$ will be either {$J=z^2$ representing altitude or altitude-related signals}, or \textcolor{black}{$J = \dot{w}^2$ representing acceleration or acceleration related signals} as shown in Figure \ref{fig:ESC_customized_structure}. This sets up the ES system to steer the insect towards the minimum of $J$ by updating the flapping angle based on updating $\tau$. 
 % It is important to note that in classic ES systems \cite{krstic2000stability,ariyur2003real} high and/or low pass filters can be used optionally. In our natural hovering ES system (Figure \ref{fig:ESC_customized_structure}), no filters at all are used. 
 Now, we are in a position to provide the proposed natural hovering ES system as 
% We propose the modeling and control system by the classic ES system for the hovering flapping flight as shown in figure \ref{fig:ESC_customized_structure}. Equivalently, the hovering flapping flight mathematical model and control system is provides 
in (\ref{eq:eq1_customized}) and (\ref{eq:eq2_customized}):  
\begin{equation}
\bm{\dot{x}} = f(\bm{x}, \tau),
\label{eq:eq1_customized}
\end{equation}
\begin{equation}
\dot{\hat{\tau}} = K J a \, \Omega \cos(\Omega t),
\label{eq:eq2_customized}
\end{equation}
where \(\bm{\dot{x}} = f(\bm{x}, \tau)\) is the insect flight dynamics model as in (\ref{eq:Model_state_space}), $K$ is the integrator gain (learning rate from feedback), $a$ is the amplitude of the modulation/demodumation signals $\Omega cos(\Omega t)$, and $\Omega$ is the natural flapping frequency of the insect/hummingbird. The natural hovering ES system can be written in state space representation as in (\ref{eq:Full_model_ESC}):
\begin{equation}
\label{eq:Full_model_ESC}
\frac{d}{dt}
\begin{pmatrix}
z \\
\phi\\
w \\
\dot{\phi} \\
\hat{\tau}
\end{pmatrix}
=
\begin{pmatrix}
w \\
\dot{\phi}\\
g - k_{d1} \, \left|\dot{\phi}\right| \, w - k_L \dot{\phi}^2 \\
-k_{d2} \, \left|\dot{\phi}\right| \, \dot{\phi} - k_{d3} \, w \dot{\phi} + \frac{\hat{\tau}}{I_F} + \frac{a}{I_F} \, \Omega \cos(\Omega t) \\
J K \frac{a}{I_F} \, \Omega \cos(\Omega t)
\end{pmatrix}
\end{equation}

\subsection{Results and Simulation}
In this Subsection, we present simulations for the proposed natural hovering ES system in \eqref{eq:Full_model_ESC} (see also Figure \ref{fig:ESC_customized_structure}) for five different flapping insects \textcolor{black}{based on literature data}, namely, hawkmoth, cranefly, bumblebee, dragonfly, hoverfly, and a hummingbird. The aim of these simulations is to validate and verify the effectiveness of the proposed natural hovering ES in achieving equilibrium hovering stability, in real-time, by adjusting the flapping angle rate based on feedback of sensations (objective function measurements).

\begin{figure*}[ht]
    \centering
    \includegraphics[width=1\linewidth]{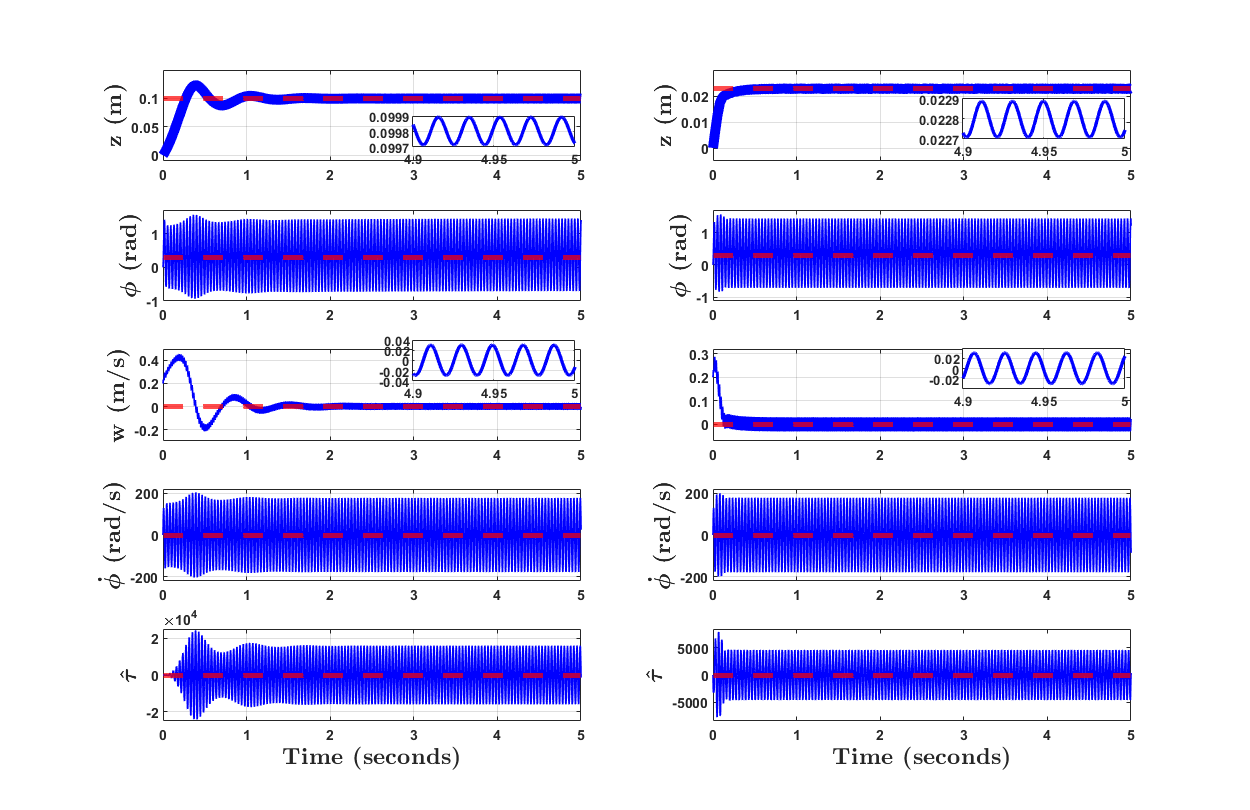}
    \caption{\textcolor{black}{State variables response of the natural hovering ES system \eqref{eq:Full_model_ESC} for \textit{\textbf{hawkmoth}}. The left column represents results under the objective function \( J = z^2 \), while the right column represents results under the objective function \( J = {\dot{w}}^2 \). The ES system (in black) stabilizes by oscillating about the hovering equilibrium (red).}}
    \label{fig:Hawkmooth_states}
\end{figure*}

\begin{figure*}[ht]
    \centering
\includegraphics[width=1\linewidth]{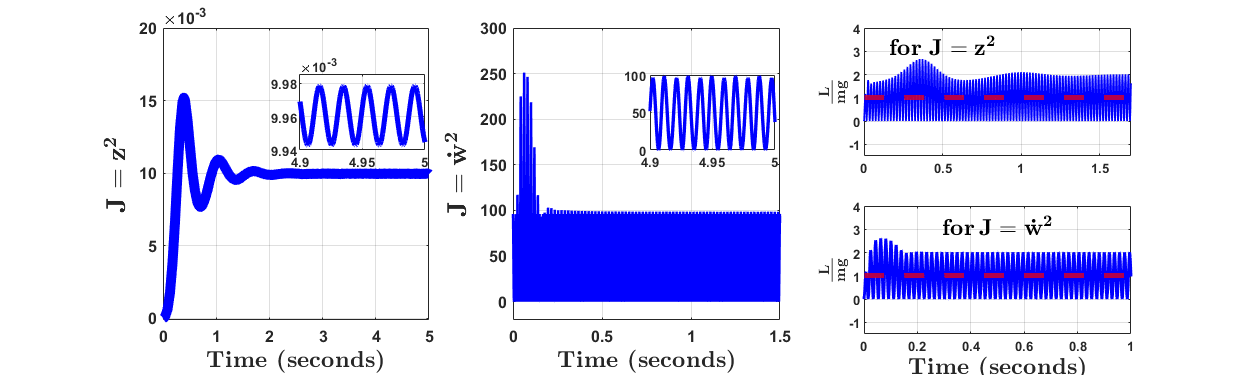}
    \caption{\textcolor{black}{For the \textit{\textbf{hawkmoth}} case, hovering condition by stabilization about a constant altitude $z$ is observed in the most left figure. Similarly, the hovering condition of balance between lift and weight is observed in the most right figure and the relevant objective function in the middle figure.}}
    \label{fig:Hawkmooth_J}
\end{figure*}

% Our simulations are conducted using two different objective functions to thoroughly evaluate the system's ability to converge to the stable hovering equilibrium, ensuring robustness and effectiveness across varying conditions.
% In this subsection, we present simulations for the proposed extremum seeking control (ESC) customized for hovering flight in flapping insects and hummingbirds. The aim of these simulations is to validate and verify the effectiveness of the proposed extremum seeking control (ESC) as a real-time feedback controller for achieving equilibrium hovering stability. Specifically, the ESC system is assessed in terms of its ability to maintain stable hovering by utilizing the natural lift generated by wing flapping. The simulations were conducted using two different objective functions to thoroughly evaluate the system's ability to converge to the stable hovering equilibrium, ensuring robustness and effectiveness across varying conditions.
In equilibrium hovering, the velocity $w$ has to be 0 as this translates into constant height/altitude $z$ (see the first differential equation in \eqref{eq:Full_model_ESC}). Therefore, in our simulations we start with an initial disturbance in $w$, that is $w_0 = 0.2$. This disturbance challenges the natural hovering ES system as it needs to overcome temporary deviations from equilibrium and re-establish hovering equilibrium (i.e., $w=0$), demonstrating stability in real-time. The simulations were performed using two objective functions: \( J = z^2 \), where its minimum represents the hovering condition (constant $z$) and \textcolor{black}{\( J = \dot{w}^2 \)}, where its minimum represents the hovering condition \textcolor{black}{(lift to weight ratio around one).}
% We started the simulation with an initial disturbance in \( w = 0.2 \) to simulate realistic flight dynamics and to model a perturbation that may occur in practical flight conditions. This disturbance challenges the system to overcome temporary deviations from stable hovering and re-establish hovering equilibrium. To assess the system's response under different control strategies, the simulations were performed using two objective functions: \( J = z^2 \), where \( z \) is the altitude, aiming to maintain altitude stability during hovering, and \( J = \left(\frac{L}{mg} - 1\right)^2 \), where \( L \) represents lift and \( mg \) is the gravitational force, focusing on stabilizing the lift-to-weight ratio around one, the equilibrium condition for hovering.
By employing both objective functions in our simulations, the system demonstrated an impressive ability to re-stabilize itself, with the vertical velocity state \( w \) oscillating around zero and the altitude remaining constant over time. This behavior reflects stable return to equilibrium hovering after disturbance. The state variables response of the natural hovering ES system \eqref{eq:Full_model_ESC} are provided in Figure \ref{fig:Hawkmooth_states} for the hawkmoth, and Figures \ref{fig:Cranefly_states}, \ref{fig:Bumblebee_states}, \ref{fig:Dragonfly_states}, \ref{fig:Hoverfly_states} and \ref{fig:Hummingbird_states} in the Appendix for the other insects and the hummingbird. The left column of each figure corresponds to results under the objective function \( J = z^2 \), while the right column of each figure corresponds to results under the objective function \textcolor{black}{\( J = \dot{w}^2 \)}. The stable hovering conditions of our simulation results (constant $z$ and lift to weight ratio of one) are further clarified in Figure \ref{fig:Hawkmooth_J} for hawkmoth and Figures \ref{fig:Cranefly_J}, \ref{fig:Bumblebee_J}, \ref{fig:Dragonfly_J}, \ref{fig:Hoverfly_J} and \ref{fig:Hummingbird_J} in the Appendix for the other insects and the hummingbird.

% In our simulations results, we modeled the lift generated by wing flapping during hovering as \( L = m K_L \dot{\phi}^2 \), where \( m \) represents the total mass of each insect, \({\phi} \) is the flapping angle of the wing, and \( K_L \) is the wing lift coefficient specific to each insect.

\textcolor{black}{Throughout the simulation results, \( \Omega = 2\pi f \). The learning-rate gain $K$ used for all simulations is given in Table \ref{Tab: ESC_model_parameters}, and the values for each insect's flapping frequency are provided in Table \ref{tab:morphological_parameters}. It is important to highlight that we used the same learning-rate gain \(K\) for all simulations across all species. This underscores the fact that the natural hovering ES is robust and requires no gain re-tuning when experiencing large variations (e.g., an individual insect/bird that undergoes wing damage to
retain its ability to hover stably despite the wing’s altered force-kinematic relationships).} Other parameters necessary to implement our natural hovering ES is given in Table \ref{Tab: ESC_model_parameters}.

\begin{table}[h!]
    \centering
    \large
\caption{\textcolor{black}{Natural hovering ES parameters for the five insects and hummingbird for different objective functions used in our simulation.}}
    \label{Tab: ESC_model_parameters}
    % \arrayrulecolor{black} % Set table border color to black
    \renewcommand{\arraystretch}{1.2} % Adjust row height for better readability
    \setlength{\tabcolsep}{8pt} % Adjust column spacing for better fit
    \begin{tabular}{|c|cc|cc|}
        \hline
        % \rowcolor{black!25} % Light black header row
        \multirow{2}{*}{\textbf{\textcolor{black}{Insect}}} & \multicolumn{2}{c|}{\textbf{\textcolor{black}{for $J = z^2$}}} & \multicolumn{2}{c|}{\textbf{\textcolor{black}{for $J = \dot{w}^2$}}} \\
        % \rowcolor{black!25} % Light black sub-header row
        & \multicolumn{1}{c}{\textbf{\textcolor{black}{$a$}}} & \multicolumn{1}{c|}{\textbf{\textcolor{black}{$K$}}} & \multicolumn{1}{c}{\textbf{\textcolor{black}{$a$}}} & \multicolumn{1}{c|}{\textbf{\textcolor{black}{$K$}}} \\
        \hline
        \textcolor{black}{Hawkmoth}     & \textcolor{black}{180}  & \textcolor{black}{10,000}  & \textcolor{black}{183.92}  & \textcolor{black}{0.5}  \\ \hline
        \textcolor{black}{Cranefly}     & \textcolor{black}{274}  & \textcolor{black}{10,000}  & \textcolor{black}{287.86}  & \textcolor{black}{0.5}  \\ \hline
        \textcolor{black}{Bumblebee}    & \textcolor{black}{966}  & \textcolor{black}{10,000}  & \textcolor{black}{969.7}  & \textcolor{black}{0.5}  \\ \hline
        \textcolor{black}{Dragonfly}    & \textcolor{black}{876}  & \textcolor{black}{10,000}  & \textcolor{black}{880.49}  & \textcolor{black}{0.5}  \\ \hline
        \textcolor{black}{Hoverfly}     & \textcolor{black}{737}  & \textcolor{black}{10,000}  & \textcolor{black}{740.24}  & \textcolor{black}{0.5}  \\ \hline
        \textcolor{black}{Hummingbird}  & \textcolor{black}{390}  & \textcolor{black}{10,000}  & \textcolor{black}{407.82}  & \textcolor{black}{0.5}  \\ \hline
    \end{tabular}
\end{table}
% The ESC parameters used in these simulations are detailed in table \ref{Tab: ESC_model_parameters}, ensuring consistency in the control strategy across different insects and objective functions. These results highlight the robustness and adaptability of the proposed ESC system across varying flight dynamics.

It is worth mentioning that we validated our natural hovering ES system further by conducting an additional set of simulations starting with negative disturbance \( w_0 = -0.2 \). We do not provide these results for space considerations, but we report that the natural hovering ES was able stabilize effectively, similar to the reported results earlier with $w_0 = 0.2$.
Overall, simulation results suggest that the proposed natural hovering ES system is a powerful and robust method for achieving real-time, stable, and model-free hovering flight. This align with biological observations, where feedback mechanisms enable flapping insects to stabilize their hovering flight in real-time in a stable manner. 
% as will be discussed in detail in the following subsection. 
\subsection{Re-examining Some Literature Results from the Eyes of Natural Hovering Extremum Seeking}
In this Subsection, we re-examine two important topics from the literature and show how they can be re-interpreted based on characterizing hovering as an ES system. The first topic is open-loop stability, and the second topic is flapping amplitude. 

\begin{figure}[htbp]
    \centering
\includegraphics[width=1\linewidth]{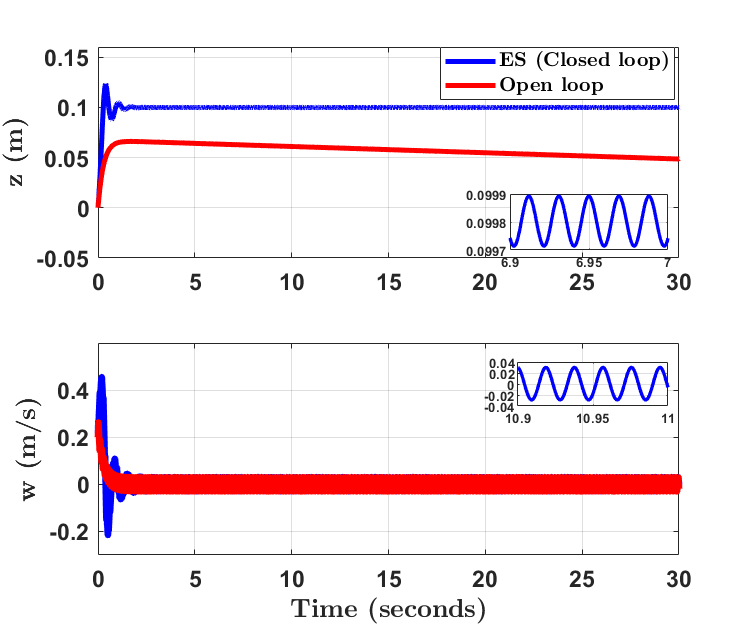}
    \caption{\textcolor{black}{Open-loop simulation (red) vs. closed-loop simulation (blue). Open-loop velocity trajectory is stabilizes about very small non-zero velocity equilibrium due to numerical errors, leading to slow linear descending in altitude. This is not the case with ES closed-loop feedback, which achieves successfully, and able to maintain, stable hovering.}}
\label{fig:open_loop_comparison}
\end{figure}

The first topic: re-interpretation of open-loop stability in hovering. \textcolor{black}{It is widely believed in the literature that hovering of flapping insects (or hummingbirds) is inherently unstable in the open-loop sense due to the pitching motion. However, recently, the authors in \cite{taha2015need} argued that open-loop stability can be shown in flapping insects while including the pitching motion for the hawkmoth. They also acknowledged \cite{taha2020vibrational} the plausibility of stability being achieved due to some feedback mechanisms similar to what was asserted in many other works in literature (see relevant discussion and citations in Subsection III.C).
While we anticipate that natural hovering ES can potentially lead to stable hovering even when pitching motion is included, we focus our attention here on open-loop stability as presented in \cite{hassan2018combinedguidancejournal} since we borrowed our model from that work.  
% However, they did so only for the hawkmoth. 
It is important to note here that there are two kinds of stability that can be achieved in open-loop setting \cite{hassan2018combinedguidancejournal} (similar observations also hold even when pitching motion \cite{taha2015need} is included). First, open-loop stability for $w=\text{constant}$. This kind of stability does not lead to hovering balance since $w=\text{constant}$ translates into linearly ascending/descending $z$ as noted in multiple figures in \cite{hassan2018combinedguidancejournal,taha2015need}. Second, open-loop stability for $w=0$, which leads to hovering balance since $w=0$ translates into constant altitude, i.e., $z=constant$ as expected in stable hovering. As a result, it is quite challenging to achieve precise stable hovering in an open-loop sense because one has to find a trim condition characterizing a flapping amplitude that leads to exact stable hovering condition (i.e., $w=0$), which can be hard to find and requires advanced methods of numerical treatments/approaches as done in \cite{hassan2018combinedguidancejournal}. Any slight deviation from said trim condition, will lead to imprecise hovering in the altitude $z$.}   
% stabnility note here, is that in many of the simulation results in \cite{taha2015need}, open-loop stability was achieved} \textcolor{black}{for the pitching motion}, but not for a hovering condition in the sense that $w=0$; open-loop stability was achieved for $w=\text{constant}$, which translates into linearly ascending/descending $z$. 
On the contrary, our natural hovering ES even for the hawkmoth case (see Figure \ref{fig:Hawkmooth_states}), stabilizes about $w=0$ even under disturbances, which always achieves a constant $z$ as naturally observed in hovering. The disparity between some of the results in  \textcolor{black}{\cite[Figures 3,4,5,6]{hassan2018combinedguidancejournal}} and our results (see Figure \ref{fig:Hawkmooth_states}) is due to the breaking (de-activation) of the sensation feedback; this leads to a system -- even if stabilized in the open-loop sense -- that is \textcolor{black}{challenging} to \textcolor{black}{maintain} stability about the required hovering condition. This is the opposite of the natural hovering ES results. Although this should be obvious, we note that the open-loop approach requires the use of kinematic equation representing flapping; if we set $K=0$ (see Figure \ref{fig:ESC_customized_structure}) in the natural hovering ES, we get an open-loop system as the modulation signal is then the kinematic equation itself that represent flapping.   
% In this subsection, we delve deeply into the validation and verification of our results by comparing them with previous work found in the literature. We also examine how (ES) system aligns with natural observations. Furthermore, we provide an in-depth exploration of the definition of hovering equilibrium stability, particularly in the context of our simulation results. Lastly, we discuss the characteristics of (ES) systems, highlighting their effectiveness as demonstrated through the simulation results.
% To begin validating our ESC results, we first compare them with the hovering stability case discussed in \cite{taha2015need}, along with other studies such as \cite{sun2014insect,taylor2003dynamic,taylor2005nonlinear}, which concluded that open-loop flapping hovering for insects, particularly hawkmoths, is inherently unstable. Although \cite{taha2015need} demonstrated open-loop stability for the hawkmoth using higher-order averaging, the system ultimately deviated from equilibrium hovering under real-world conditions, as shown in their simulation results.
To clarify further our re-interpretation of the open-loop stability, 
% results in \textcolor{red}{\cite{taha2015need}} \textcolor{red}{\cite{hassan2018combinedguidancejournal}}
we conducted a simulation comparing between the open-loop case \textcolor{black}{from \cite{hassan2018combinedguidancejournal} by using exact $U_{trim}$ value (amplitude of flapping) computed numerically from optimized shooting method} and our natural hovering ES 
% de-activated feedback (i.e., $K=0$ in Figure \ref{fig:ESC_customized_structure}) 
with closed-loop ES feedback system (i.e., $K>0$). We used the objective function \( J = z^2 \) for the natural hovering ES (closed-loop).
% , starting with a initial disturbance in the vertical velocity, \( w_0 = 0 \, \text{m/s} \). 
The simulation results (see Figure \ref{fig:open_loop_comparison}) shows that the trajectory of the ES closed-loop feedback system diverges away from its alignment with the open-loop trajectory, seeking stabilization that achieves hovering (oscillating about $w=0$ and constant $z$). However, the open-loop trajectory \textcolor{black}{is struggling to maintain $z=constant$ even with very small numerical errors in the computation of hovering balance flapping amplitude. With lack of feedback, this slight divergence of $z$ cannot be corrected.} 
% This comparison validates the effectiveness of the ES closed-loop feedback system in steering the system back to the equilibrium hovering condition, even in the face of large initial disturbances. This behavior is more representative of the real-world conditions faced by insects and hummingbirds, making the ESC system a more reliable and realistic control strategy.

The second topic: flapping angle amplitude. We re-examined how the flapping angle amplitude of our natural hovering ES compares to natural experimental observations. That is, we check how $\phi$ oscillates under the natural hovering ES system, and compare the amplitude of its oscillation \( \Phi^\circ \) with reported values in the litrature based on natural/experimental observation \cite{145,ellington1984aerodynamics,dudley1990mechanics}. We provide said comparison in Table \ref{tab:comparison of flapping amplitude values}. 
 % as shown in table \ref{tab:morphological_parameters}.
The comparison reveals only slight differences between the natural observations and the flapping angle amplitudes produced by the natural hovering ES, demonstrating that the real-time natural hovering ES system is effective in mimicking the natural hovering flight behavior.
\begin{table}[htbp]
    \centering
    \large
    \caption{Flapping angle amplitude values reported from natural/experimental observations (middle column) vs. our natural hovering ES system for the five insects and hummingbird used in this work.}
    \begin{tabular}{|c|c|c|}
        \hline
    \textbf{Insect} & \textbf{ $\bar{\Phi}^o$ (rad)} & \textbf{ES $\bar{\Phi}^o$ (rad)} \\
        \hline
        Hawkmoth & 1.05 & 1.07 \\ \hline
        Cranefly & 1.07 & 1.03 \\ \hline
        Bumblebee & 1.01 & 0.98 \\ \hline
        Dragonfly & 0.95 & 0.91 \\ \hline
        Hoverfly & 0.78 & 0.75 \\ \hline
        Hummingbird & 1.22 & 1.19 \\ \hline
    \end{tabular}
    \label{tab:comparison of flapping amplitude values}
\end{table}

% Overall, the simulation results demonstrate that the proposed ESC system effectively stabilizes the hovering equilibrium in flapping insects and hummingbirds using different objective functions. The ES system consistently brought the system back to hovering equilibrium condiiton. Moreover, the slight difference between the ESC results and natural observations further confirms the system's accuracy, making it a promising real-time feedback controller for flapping hovering insects and hummingbirds.

\subsection{\textcolor{black}{Comparison Between PID and Natural Hovering Extremum Seeking Approach}}
\textcolor{black}{
In this Subsection, we move forward with a comparison between the proposed natural hovering ES system and Proportional-Integral-Derivative (PID) approach for controlled hovering. The reason we choose a comparison with PID controller is that it was used in hovering literature (eg., \cite{Delft_Conference,DelftDissertation2014,feedback25ristroph2010discovering}), but it is also a traditional simple controller that has the potential of operating in a model-free \textcolor{black}{and real-time} fashion with low computational cost. We believe this comparison is important as it highlights for many readers, especially those not familiar with control theory/systems, the advantages and robustness of the proposed natural hovering ES system.}

\textcolor{black}{PID controllers rely on error feedback from a \textit{predefined} reference signal and are widely, and most effectively, used for linear time-invariant (LTI) systems where they can be considered model-free and have much larger tolerance in tuning their gains. However, since flapping-wing dynamics are nonlinear, time-varying systems, PID controllers cannot be applied directly or naturally on such systems. Consequently, applying PID control to flapping systems requires multiple approximations and extensions known in control theory and systems. In fact, it is a general knowledge in control theory and systems that PID control application in highly nonlinear systems is suboptimal and leads to gains that are: (1) dependent on a predefined equilibrium and linearization of the system around said equilibrium; and (2) very sensitive to tuning with sometimes very small margin of tolerance. It is worth noting that flapping-wing systems are not just highly nonlinear systems, but they are time-varying as well, which adds another layer of complexity to the problem of controlling hovering via PID-based technique.} 

\begin{figure}[ht]
    \centering
    \includegraphics[width=1\linewidth]{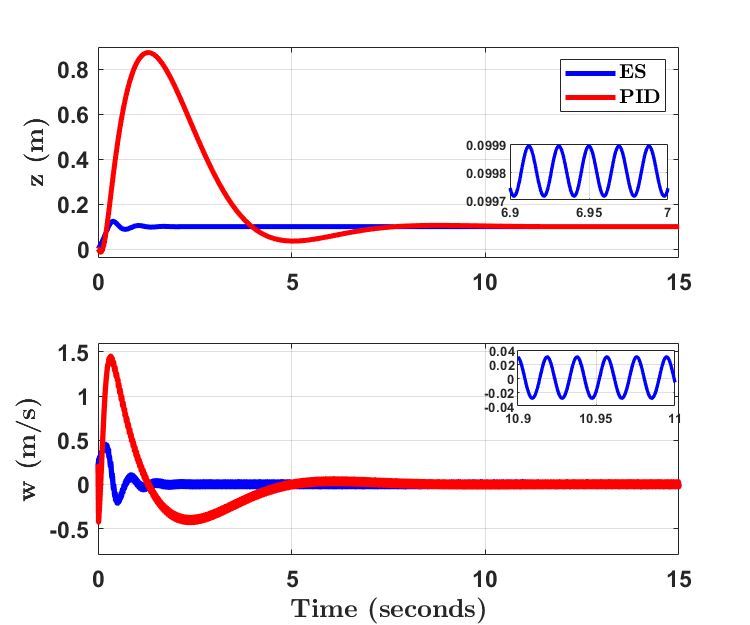}
    \caption{\textcolor{black}{Comparison of altitude stabilization using PID (red) vs. natural hovering ES system (blue). Natural hovering ES system is performing better in terms of settling time and its lesser overshoot.}}
    \label{fig:pid}
\end{figure}

\textcolor{black}{Next, we construct a proper PID control for hovering. In our system dynamics (\ref{eq:Model_state_space}), the flapping amplitude \( a \) is used as the control input for the PID controller (inline with literature), while the flapping torque is expressed as \( \tau = a \cos(\omega t) \). The PID controller operates by minimizing the error feedback signal, which is defined as \( e = z - z_{\text{desired}} \), ensuring that the system stabilizes at the desired altitude.}
\textcolor{black}{To construct the PID design and find its gains, the following steps are taken \cite{pid_book}:
\begin{enumerate}
    \item We convert the system from being nonlinear time-varying to nonlinear time-invariant. This can be achieved via variation-of-constant averaging methodology in a similar fashion to \cite{hassan2018combinedguidancejournal}.
    \item We linearize the resulting nonlinear time-invariant system from step (1) to obtain an LTI system.
    \item We derive an open-loop transfer function between the input \( a \) and the output \( z \).
    \item We follow standard PID tuning techniques as detailed in \cite{pid_book}.
\end{enumerate}}
% the system must first be converted into a linear time-invariant (LTI) form. However, our flapping flight model is inherently a nonlinear time-periodic system, meaning that:
% 1. The first step requires obtaining the averaged system(Nonlinear Time-Varying system), following the methodology in \cite{hassan2018combinedguidancejournal}.
% 2. The second step involves linearizing this averaged system to derive an open-loop transfer function between the input \( a \) and the output \( z \).
% 3. Following standard PID tuning techniques from \cite{pid_book},
We obtained the PID gains as error proportion gain \( k_p = 3.66e-03 \), integral of error gain \( k_i = 2.936e-03 \), and derivative of error gain \( k_d = 2.00e-03 \). We would like to report that the small derivative of error gain $k_d$ was added to reduce overshoot even though it did not come out naturally from the tuning process. \textcolor{black}{We also would like to report that the system is very sensitive to the gains. To illustrate this sensitivity, the initial disturbance \(w_0\) was varied while the PID gains (mentioned earlier) kept constant.
Table~\ref{tab:ES_PID_disturbance} compares overshoot and settling time for the resulting PID response against the natural hovering ES controller, which uses the same learning-rate gain \(K\) in every case. Even small changes in \(w_0\) drive the PID overshoot above \(7\!\times\!10^{2}\,\%\) and lengthen the settling time to \(\approx10\) s, whereas ES maintains overshoot below \(26\%\) and settles in time \(\le1.6\) s without re-tuning. It should be noted that the settling time and overshoot values were calculated using MATLAB’s built-in (stepinfo) function. To visualize this effect, Figure \ref{fig:pid} shows the response of the PID-controlled hovering versus the natural hovering ES system for $w_0=0.2$, which shows that the natural hovering ES system is performing better in terms of settling time and its lesser overshoot. These results confirm that the PID controller must be re-optimized for each disturbance magnitude, while the ES law is robust to large variations/disturbances. For the full details of the PID tuning methodology and derivations, we provide the implementation files in \cite{githubdae}}.
\begin{table}[htbp]
    \centering
    \large
    \caption{\textcolor{black}{Overshoot and settling‐time performance of the ES and PID controllers for five initial perturbations \(w_0\).}}
    \label{tab:ES_PID_disturbance}
    \textcolor{black}{
    \begin{tabular}{|c|c|c|c|c|}
        \hline
        \multirow{2}{*}{\(\boldsymbol{w_0}\)} &
        \multicolumn{2}{c|}{\textbf{Overshoot [\%]}} &
        \multicolumn{2}{c|}{\textbf{Settling time (s)}} \\ \cline{2-5}
        & ES & PID & ES & PID \\ \hline
        0.2 & 23.63 & 770.22 & 1.39 & 10.22 \\ \hline
        0.3 & 25.22 & 769.85 & 1.39 & 10.20 \\ \hline
        -0.2 & 22.82 & 772.45 & 1.54 & 10.26 \\ \hline
        -0.3 & 23.67 & 773.28 & 1.60 & 10.27 \\ \hline
        -0.1 & 22.41 & 771.70 & 1.29 & 10.24 \\ \hline
    \end{tabular}}
\end{table}
\textcolor{black}{For more demonstrative, stronger, and fairer comparison, we tried to get the optimal PID gains by using the genetic algorithm (GA) optimizer. Running GA for each insect, the optimal PID gains are obtained as shown in Table~\ref{tab:PID_gains}; the GA code and configuration files are available in~\cite{githubdae}. To demonstrate the sensitivity of even optimal PID gains to changes in operating conditions, the GA was rerun for the hawkmoth under five different initial \(w_0\).  
Table~\ref{tab:Hawkmoth_PID_tuning} shows that every change in \(w_0\) yields different \((k_p,k_i,k_d)\) values, confirming that the PID controller must be re-optimized (re-tuned) whenever the operating condition changes.}

\begin{table}[htbp]
    \centering
    % \normalsize
    \large
    \caption{\textcolor{black}{Optimal PID gains obtained by GA for $w_0=0.2$}}
    \label{tab:PID_gains}
    \textcolor{black}{
    \begin{tabular}{|c|c|c|c|}
        \hline
        \textbf{Insect} & \(k_p\) & \(k_i\) & \(k_d\) \\ \hline
        Hawkmoth    & \(8.912\text{e}{-02}\) & \(2.585\text{e}{-01}\) & \(2.405\text{e}{-02}\) \\ \hline
        Bumblebee   & \(5.000\text{e}{-04}\) & \(4.000\text{e}{-04}\) & \(1.000\text{e}{-04}\) \\ \hline
        Cranefly    & \(1.421\text{e}{-04}\) & \(3.983\text{e}{-04}\) & \(4.679\text{e}{-05}\) \\ \hline
        Dragonfly   & \(5.000\text{e}{-04}\) & \(4.000\text{e}{-04}\) & \(9.999\text{e}{-05}\) \\ \hline
        Hoverfly    & \(3.337\text{e}{-04}\) & \(4.000\text{e}{-04}\) & \(9.999\text{e}{-05}\) \\ \hline
        Hummingbird & \(5.000\text{e}{-02}\) & \(4.000\text{e}{-02}\) & \(1.000\text{e}{-02}\) \\ \hline
    \end{tabular}}
\end{table}

\begin{figure}[ht]
    \centering
    \includegraphics[width=\linewidth]{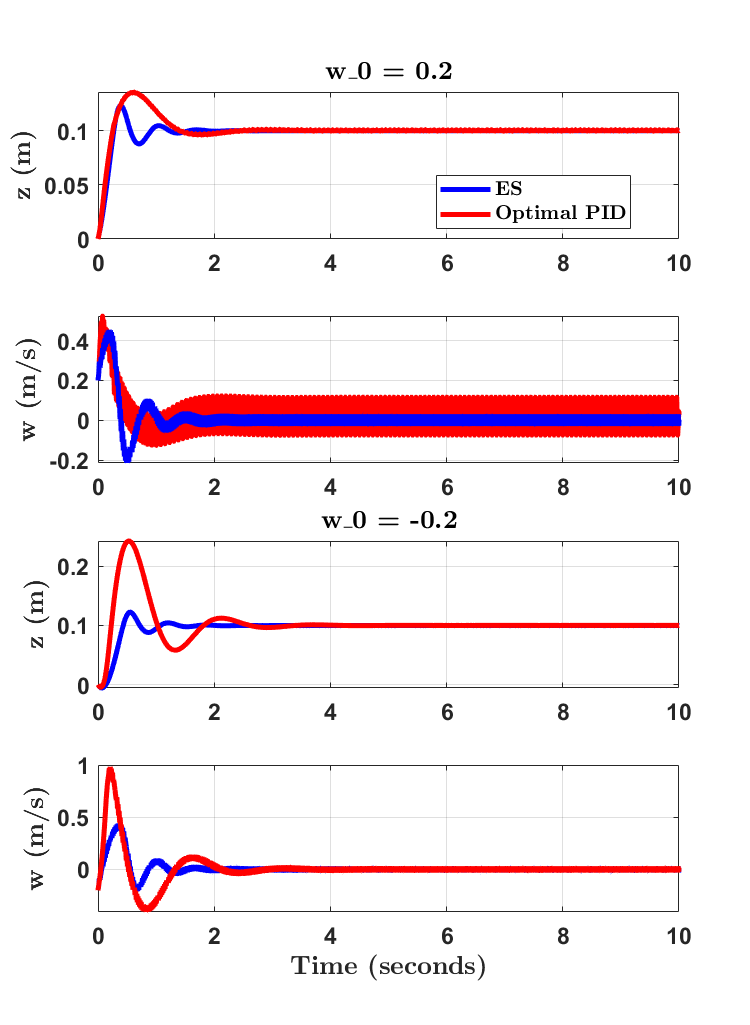}
    \caption{\textcolor{black}{Comparison of altitude stabilization using the optimal PID controller (red) versus the natural hovering ES system (blue).
             The ES loop settles faster and with smaller overshoot for both initial conditions (\(w_{0}=0.2\) and \(w_{0}=-0.2\)), while using the same ESC feedback gain in every case and different optimal PID gains as shown in table \ref{tab:Hawkmoth_PID_tuning} .}}
    \label{fig:comparison_ES_optimal_PID}
\end{figure}

\begin{table}[htbp]
    \centering
   \large
    \caption{\textcolor{black}{Hawkmoth—optimal PID gains identified for five initial disturbances \(w_0\).}}
    \label{tab:Hawkmoth_PID_tuning}
    \textcolor{black}{
    \begin{tabular}{|c|c|c|c|}
        \hline
        \(\boldsymbol{w_0}\) & \(\boldsymbol{k_p}\) & \(\boldsymbol{k_i}\) & \(\boldsymbol{k_d}\) \\ \hline
        \(0.20\) & \(8.912\mathrm{e}{-02}\) & \(2.585\mathrm{e}{-01}\) & \(2.405\mathrm{e}{-02}\) \\ \hline
        \(0.30\) & \(1.147\mathrm{e}{-01}\) & \(3.178\mathrm{e}{-01}\) & \(2.429\mathrm{e}{-02}\) \\ \hline
        \(-0.10\) & \(2.411\mathrm{e}{-02}\) & \(1.027\mathrm{e}{-01}\) & \(2.296\mathrm{e}{-02}\) \\ \hline
        \(-0.20\) & \(2.386\mathrm{e}{-02}\) & \(1.152\mathrm{e}{-01}\) & \(6.272\mathrm{e}{-03}\) \\ \hline
        \(-0.30\) & \(2.188\mathrm{e}{-02}\) & \(8.514\mathrm{e}{-02}\) & \(4.248\mathrm{e}{-03}\) \\ \hline
        \end{tabular}}
\end{table}
\textcolor{black}{Moreover, Table~\ref{tab:ES_optimal_PID_disturbance} compares overshoot and settling time for the natural hovering ES (same gain \(K\) in all cases) and the GA–tuned (i.e., optimal) PID controller. The optimal PID improves on the hand-tuned version when compared to its performance in Table \ref{tab:ES_PID_disturbance}, yet still requires more time to settle with high overshoot. In contrast, the natural hovering ES remains more efficient and robust with less overshoot and settling time without any re-tuning. To visualize this effect, Figure~\ref{fig:comparison_ES_optimal_PID} plots altitude \(z(t)\) and velocity \(w(t)\) for two representative disturbances, \(w_0=0.2\) and \(w_0=-0.2\).  
Even with their GA-optimized gains (shown in Table(\ref{tab:ES_optimal_PID_disturbance})), the optimal PID response takes longer to converge and overshoots more than the natural hovering ES, which further confirms the same conclusion about the ES effectiveness and robustness. These results highlight a core advantage of the natural hovering ES controller: a single gain setting delivers consistent performance across a wide range of initial conditions/disturbances, whereas PID requires case-by-case re-optimization with sub-optimal performance/robustness. We provide the implementation files in \cite{githubdae}}

\begin{table}[htbp]
    \centering
    \large
    \caption{\textcolor{black}{Overshoot and settling‐time performance of the ES and optimal PID controllers for five initial perturbations \(w_0\).}}
    \label{tab:ES_optimal_PID_disturbance}
    \textcolor{black}{
    \begin{tabular}{|c|c|c|c|c|}
        \hline
        \multirow{2}{*}{\(\boldsymbol{w_0}\)} &
        \multicolumn{2}{c|}{\textbf{Overshoot [\%]}} &
        \multicolumn{2}{c|}{\textbf{Settling time (s)}} \\ \cline{2-5}
        & ES & Optimal PID & ES & Optimal PID \\ \hline
        0.2 & 23.63 & 34.94 & 1.39 & 2.27 \\ \hline
        0.3 & 25.22 & 29.98 & 1.39 & 1.97 \\ \hline
        -0.2 & 22.82 & 142.88 & 1.54 & 3.18 \\ \hline
        -0.3 & 23.67 & 178.51 & 1.60 & 2.45 \\ \hline
        -0.1 & 22.41 & 68.36 & 1.29 & 8.85 \\ \hline
    \end{tabular}}
\end{table}

Overall, the key advantages of the natural hovering ES system over PID control stem from its model-free, real-time adaptability. Unlike PID when it is applied to nonlinear systems, which requires a linearized transfer function, ES does not rely on explicit system knowledge by any mean. Additionally, as discussed earlier, natural hovering ES does not require an additional input to control the system, but the flapping oscillation itself is the control input. 
Additionally, PID control struggles with unknown objectives (even if it is able to handle model-free dynamics), as it requires a predefined reference altitude \( z_{\text{desired}} \), whereas ES autonomously seeks the optimal hovering condition corresponding to a maximum or minimum of an unknown objective function (refer to earlier discussions about light source seeking). That is, ES is capable of performing real-time optimization and source seeking, unlike PID, which relies purely on error minimization. lastly, PID controllers require precise tuning that may need to be re-adjusted when controlling nonlinear systems, especially if the controller is stabilizing about different equilibrium (altitudes) from a time to another, whereas ES naturally adapts in real-time without requiring manual gain adjustments.
% \textcolor{black}{It should be noted also that the lack of need for gain tuning is an important strength of ES—not across species, but as a potential way for an individual animal that, for example, undergoes wing damage to retain its ability to hover stably despite the wing’s altered force–kinematic relationships.} 
% \textcolor{black}{we should emphasize also that for the PID controller we have to tune 3 gains basically which needs some many and complicated steps that depend on the averaging system (in high dimension) as explained in \cite{pid_book}, or they will need offline tunning to get the optinal values as we did using GA or as done before using neural network in \cite{neural}. However, our proposed ES needs mainly 3 parameters to tune, frequency f, amplitutde a, and learning rate k gain. However, due to the natural flapping problem and based on the bilogical observations, the f frequency is determined by the natural flapping frequency of each insect and same also for a that detrmined by the natural flapping ampltitude values (see Table IV for comparison), the only remaining value can should be tune real time without going into complicated steps, is K which is also is a constant value for all insects as showin in table III }
\textcolor{black}{To summarize the fundamental differences between PID control and the natural hovering ES system, we present a comparative analysis in Table \ref{table:pid_vs_es}.}

\begin{table*}[ht]
\centering
\renewcommand{\arraystretch}{1.3} % Adjust row height for better readability
\setlength{\tabcolsep}{3pt} % Reduce column spacing
\caption{\textcolor{black}{Comparison between PID hovering control and natural hovering ES system}}
\label{table:pid_vs_es}
 % Set table border color to black
\begin{tabular}{|p{4cm}|p{6cm}|p{6cm}|} 
\hline
\textcolor{black}{Feature} & \textcolor{black}{PID Control} & \textcolor{black}{Natural Hovering ES System} \\ 
\hline
\textcolor{black}{Knowledge of the dynamic system model} & \textcolor{black}{May require an explicit/partial knowledge for the dynamic system model if it is nonlinear, at least for re-tuning.} & \textcolor{black}{Entirely model-free with fixed leaning-rate gain.} \\ 
\hline
\textcolor{black}{Tuning process} & \textcolor{black}{Needs system linearization  
and parameter tuning that may need re-adjustments if the controller stabilizes about different equilibrium from a time to another.} & \textcolor{black}{Self-adaptive in real-time and its tuning is operable for wide range of conditions.} \\ 
\hline
\textcolor{black}{Handling of unknown objectives} & \textcolor{black}{Requires access to predefined  
reference signal \( z_{\text{desired}}.\). That is, it needs global information.} & \textcolor{black}{It does not require explicit knowledge of the objective function (including its maximum or minimum information). It operates with local measurements} \\ 
\hline
\textcolor{black}{Computational complexity} & \textcolor{black}{For linear systems, it needs an operation of proportional feedback, an operation of differentiation of error, and an operation of integration of error. For nonlinear systems, all the aforementioned is included in addition to some other operations based on the model (e.g., re-adjusting tuning).} & \textcolor{black}{Minimal computation as it needs only one integrator operation.} \\ 
\hline
\end{tabular}
\end{table*}

\subsection{\textcolor{black}{Effect of Delay and Disturbance on the Natural Hovering ES System}}
\textcolor{black}{
In real-world, delays and disturbances/noises are inevitable and can significantly impact the performance of any control system. However, ES systems are usually able to accommodate stochasticity in input or output (feedback measurements), including delays -- see, for example, some recent developments \cite{oliveira2016extremum,oliveira2022extremum,ruvsiti2018stochastic,liu2012stochastic} and some experimental results involving ES systems in noisy real-world robotic applications of light source seeking \cite{ECC2024,unicycle_bio_inspiration}. In this Subsection, we demonstrate the robustness of the natural hovering ES system under some uncertain conditions, namely the effects of time delay and noise in feedback measurements.}

\begin{figure}[ht]
    \centering
    \includegraphics[width=1\linewidth]{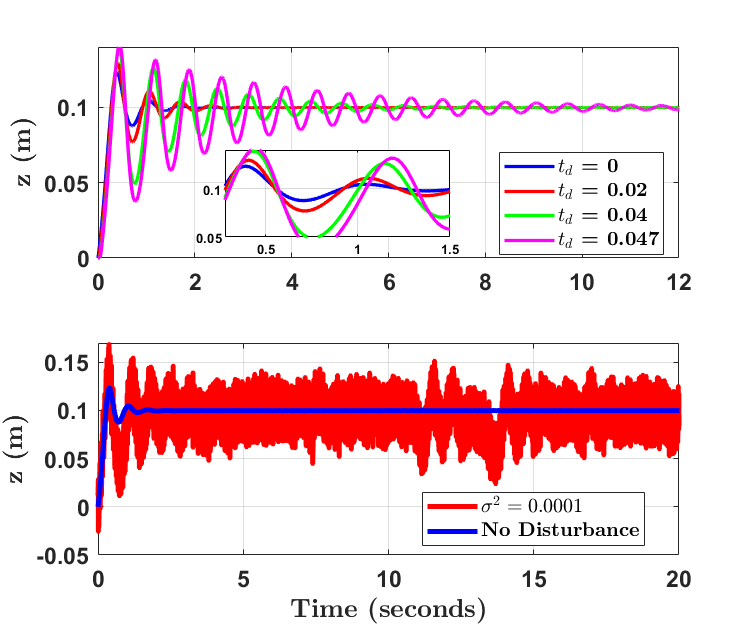}
    \caption{\textcolor{black}{Effect of delay (top) and persistence noise/disturbance (bottom) on the hovering ES system. The system remains stable under delays and persistent noise/disturbances.}}
\label{fig:delay_dis}
\end{figure}

\textcolor{black}{System delays often arise due to sensor latency, actuator response times, or computational processing delays. To analyze the impact of time delay, we introduce a fixed delay \( t_d \) in the altitude \( z \) feedback loop and examine the system’s response, in a methodological similarity to relevant processes in \cite{delay}. Consequently, the altitude \( z \) directly influences the objective function, expressed as \( J = z^2 \). The delay values tested in this study, on the hawkmoth case, are \( t_d = 0 \) (no delay), \( t_d = 0.02 \) seconds, \( t_d = 0.04 \) seconds, and \( t_d = 0.047 \) seconds. These values are significantly larger than the solver time step of 0.0001 seconds, meaning that each delay period spans hundreds of solver steps. Figure \ref{fig:delay_dis} (top) presents the altitude response under different delay conditions. As the delay increases, the system exhibits longer transient oscillations and a slower convergence rate. The zoomed-in section highlights how the phase shift increases with larger delays, impacting the system's ability to stabilize quickly. However, despite the added delay, the natural hovering ES system successfully maintains a stable hovering state.}

\textcolor{black}{Biological flight systems are constantly exposed to environmental disturbances, including wind gusts and air turbulence, among others. To evaluate the robustness of the natural hovering ES system against such uncertainties and persistent noise condition, we apply noise to the measurement (sensation) of the objective function in a similar fashion to \cite[Figure~12]{pokhrel2022novel}. we introduce a Gaussian noise disturbance added to the feedback altitude \( z \), which affects the objective function \( J = z^2 \), with a variance of \( \sigma^2 = 0.0001 \). The system response under no disturbance and with added noise is illustrated in Figure \ref{fig:delay_dis} (bottom). The natural hovering ES system remains stable even under noisy conditions, but small fluctuations are observed in the altitude due to the external perturbations. Despite the persistent noise, the system continues to hover within an acceptable altitude range compared to the case with no noise.}

\textcolor{black}{The simulation results in Figure \ref{fig:delay_dis} demonstrate the ability of the natural hovering ES system to handle both delays and disturbances in the feedback measurements, ensuring stable flight even under challenging conditions. While delays introduce transient oscillations and slower response, the system does not diverge or lose stability. Similarly, while persistent noise/disturbances cause small variations in altitude, the controller compensates effectively, maintaining a consistent hovering state. These findings, we believe, reinforce the biological plausibility of ES-based control for hovering, as it mimics the robustness of flapping insects and hummingbirds and their adaptability to real-world uncertainties while hovering.}

\section{Stability analysis of our natural hovering extremum seeking system and comparison with open-loop}

\label{sec:Stability}
Even though ES systems generally have stability guarantees and are known for their robustness \cite{scheinker2024100}, it is important to note that our natural hovering ES system needs more investigation to assert its stability. That is, our natural hovering ES system \eqref{eq:Full_model_ESC} has some difference when compared to the classic class of ES systems in \cite{krstic2000stability,ariyur2003real}. In particular, in the natural hovering ES system \eqref{eq:Full_model_ESC}, the modulation and demodulation signals have large amplitudes (i.e., scaled by $\Omega$) which is not the case in the classic ES structure -- see Figure \ref{fig:ESC_customized_structure} as opposed to the classic structure \cite[Figure 1]{krstic2000stability}. Hence, we proceed with our stability analysis to the natural hovering ES system \eqref{eq:Full_model_ESC} based on the framework provided in \cite{elgohary2025extremum}.  
The flight dynamics of hovering flapping-wing of insects have been traditionally treated as a nonlinear, time-periodic (NLTP) system \cite{taha_review}, where the equilibrium state is represented by a periodic orbit that oscillates about a fixed point equilibrium. Not surprising, averaging techniques \cite{maggia2020higher} have been the standard method in analyzing the NLTP system representing flapping flight and hovering. Averaging techniques are particularly appealing in such analysis because they are an analytical approach, which provides better physical understanding of flapping flight phenomena and the possibility of mimicking control designs. We note that our natural hovering ES system \eqref{eq:Full_model_ESC} is too, an NLTP system. However, for averaging techniques \cite{maggia2020higher} to apply, the NLTP system (even if it is in a control-affine form \cite{pokhrel2023higher}) needs to satisfy the averaging canonical form  
% \subsection{Averaging Theory}
% The Generalized Averaging Theory (GAT) is an extension of classical averaging theory, designed to approximate the solutions of nonlinear, time-periodic dynamical systems. It offers a method to analyze systems with time-periodic behavior by separating fast oscillations from slow dynamics. While classical averaging is typically used for systems with high-frequency oscillations, GAT generalizes this approach, allowing for higher-order approximations of more complex nonlinear systems. A key aspect of GAT is its use of Lie brackets, which describe the interaction of vector fields that represent the system’s dynamics. These higher-order terms improve the accuracy of the approximation and allow GAT to handle systems with multiple time scales. GAT provides a powerful tool for studying the stability and dynamics of nonlinear time-periodic systems by decoupling the fast oscillatory behavior from the slow-changing dynamics. In \cite{agravcev1979exponential}, the GAT initially introduced by
% Agrachev and Gamkrelidze, and later, Sarychev \cite{sarychev2007stability} and Vela \cite{vela2003averaging} expanded on their ideas to further generalize the classical averaging theorem. For more details about GAT, the reader is referred to \cite{maggia2020higher}. 
$\dot{\bm{x}} = \epsilon \bm{g}(\bm{x}, t)$,
where \( \epsilon \) is a small parameter, \( \bm{g} \) is smooth in \(x\) and a \( T \)-periodic vector field in \( t \). Now, does our natural hovering ES system satisfies said averaging canonical form? 

Let us first re-write the natural hovering ES system \eqref{eq:Full_model_ESC} in the form of a nonlinear control-affine system as:
\begin{equation}
\dot{\bm{x}} = \bm{Z}(\bm{x}) + \bm{G}(\bm{x}) \tau_{\phi},
\label{eq:control_affine_form}
\end{equation}
where $\bm{x} = \begin{bmatrix}
z & \phi & w & \dot{\phi} & \hat{\tau}
\end{bmatrix}^T$
is the state vector, the vector fields \(\bm{Z}(\bm{x})\) and \(\bm{Y}(\bm{x})\) are defined as follows:
\[
\bm{Z}(\bm{x}) = \begin{bmatrix}
w \\
\dot{\phi} \\
g - k_d1 \, | \dot{\phi} | \, w - k_L \, \dot{\phi}^2 \\
- k_d2 \, | \dot{\phi} | \, \dot{\phi} - k_d3 \, w \, \dot{\phi} + \frac{\hat{\tau}}{I_F} \\
0
\end{bmatrix},
\]
and
\[
\bm{G}(\bm{x}) = \begin{bmatrix}
0 & 0 & 0 & \frac{a}{I_F} & \frac{JKa}{I_F}
\end{bmatrix}^T,
\]
while $
\tau_{\phi} =  \, \Omega \cos(\Omega t)$.
% It is important to note that the states \( z \) and \( \phi \) are ignorable coordinates, meaning they do not affect on the system's dynamics. 
It is clear that the system \eqref{eq:control_affine_form} cannot be written in the averaging canonical form, deeming standard averaging techniques \cite{maggia2020higher,pokhrel2023higher}(e.g., direct/first-order or higher-order averaging) inapplicable. However, the system \eqref{eq:control_affine_form} satisfies the Variation of Constants (VOC) condition in \cite{bullo2002averaging,maggia2020higher,elgohary2025extremum}. Clearly, the NLTP system \eqref{eq:control_affine_form} (and hence our natural hovering ES system \eqref{eq:Full_model_ESC}) can be written as: 
\begin{equation}
\dot{\bm{x}} = \bm{Z}(\bm{x}) + \frac{1}{\epsilon} \bm{Y}\left(\bm{x}, \frac{t}{\epsilon}\right),
\label{eq:NLTP_voc}
\end{equation}
by taking $\epsilon =1/\Omega$ and $\bm{Y}\left(\bm{x}, \Omega t\right)=\bm{G}(\bm{x})cos(\Omega t)$.
As a result of the VOC condition \cite{bullo2002averaging,maggia2020higher,elgohary2025extremum} being satisfied, an averaged system to \eqref{eq:NLTP_voc} can be derived and used for stability analysis: 
% applying the first order averaging formula in (\ref{eq:Averaging NLTP},\ref{eq:high_order_averaging}), we obtain:
\begin{equation}
\dot{\bar{\bm{x}}} = \bm{Z}(\bar{\bm{x}}) + \frac{a^2}{4} \left[ \bm{Y}, [\bm{Y}, \bm{Z}] \right](\bar{\bm{x}}),
\label{eq:VOC_result}
\end{equation}
where \([\bm{Y}, \bm{Z}]\) is the Lie bracket between the vector fields \(\bm{Y}\) and \(\bm{Z}\), defined as
\begin{equation}
[\bm{Y}, \bm{Z}] = \frac{\partial \bm{Z}}{\partial \bm{x}} \bm{Y} - \frac{\partial \bm{Y}}{\partial \bm{x}} \bm{Z}.
\end{equation}
% where $\int_0^T \int_0^{t} \int_0^{s_1} \tau_{\phi}(s_2) \, \tau_{\phi}(s_1) \, ds_2 \, ds_1 \, dt = \frac{a^2}{4}$, and $\tau_{\phi} = a \Omega \cos(\Omega t)$ $a$ is the amplitude of the modulation/demodumation signals in (\ref{eq:Full_model_ESC}). $\Omega cos(\Omega t)$  .
% It is important to note that, due to the assumed cosine waveform for \( \tau_{\phi} \) in system \ref{eq:Full_model_ESC}, the following condition is satisfied:
% \[
% \int_0^T \int_0^t \tau_{\phi}(\sigma) \, d\sigma \, dt = 0
% \]
% \subsection{Dynamics Smooth Approximation }
A key limitation in applying VOC-based averaging via \eqref{eq:VOC_result} is that smothness in $\bm{x}$ is not satisfied in the vector field \( \bm{Z}(\bm{x}) \) due to the presence of nonsmooth function, namely $|\dot{\phi}|$ (see \eqref{eq:control_affine_form}). This lack of smoothness curtails the use of differentiation to compute the Lie bracket in \eqref{eq:VOC_result}. To address this, we introduce a smooth approximation for the absolute value function \( |\dot{\phi}| \) similar to \cite{taha2016geometriclong}.

\begin{figure}[htbp]
    \centering
\includegraphics[width=1\linewidth]{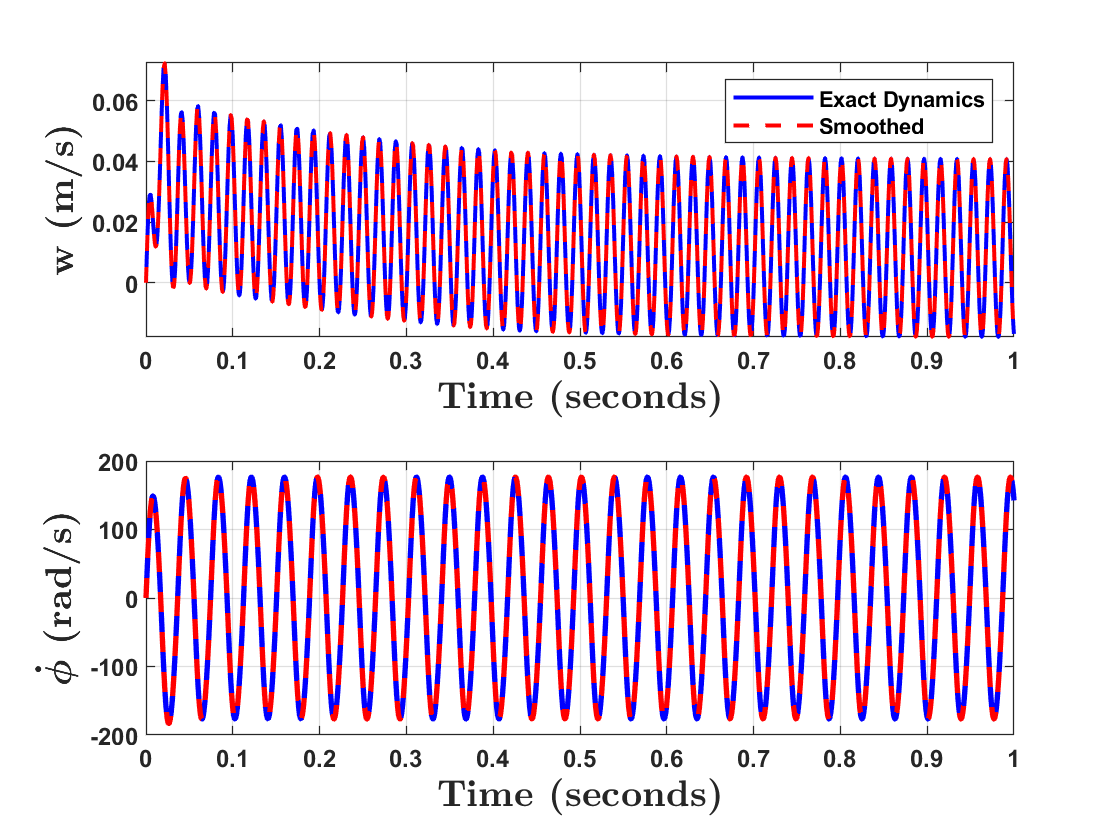}    \caption{Comparison between the responses of exact ES dynamics and approximate/smoothed dynamics.}
\label{fig:smoothing_comparison}
\end{figure}

We approximate the absolute value function as: \( |\dot{\phi}| \approx \dot{\phi} h(\dot{\phi}) \), where \( h(\dot{\phi}) = \frac{2}{\pi} \tan^{-1}(n \dot{\phi}) \). The larger the $n$, the better the smoothing approximation. Through our analysis, we found that using \( n = 50 \) works well for all insects and hummingbird, offering a sufficiently accurate approximation. Furthermore, Figure \ref{fig:smoothing_comparison} shows a precise match between the actual nonsmooth ES dynamics per \eqref{eq:Full_model_ESC} of the hawkmoth and the smoothed ES dynamics with \( n = 50 \), specifically for the states \( w \) and \( \dot{\phi} \) starting the simulation with $w_0=0$. In said simulation, we used $J= z^2$.

\begin{table*}[ht]
\centering
\large
\textcolor{black}{
\caption{\textcolor{black}{Eigenvalues of the averaged system for various insects and hummingbird under different objective functions, demonstrating stability for \( w \) and \( \dot{\phi} \).}}
\label{tab:eigenvalues}
\begin{tabular}{|c|c|c|c|c|}
\hline
Insect &
Open loop &
Optimal PID &
ESC (\(J=z^{2}\)) &
ESC (\(J=\dot w^{2}\)) \\ \hline
Hawkmoth &
\(-5.11e{+05}\) &
\(-4.81e{+05}\) &
\(-5.25e{+05}\) &
\(-5.48e{+05}\) \\
& \(-1.78e{+04}\) &
\(-1.66e{+04}\) &
\(-1.82e{+04}\) &
\(-1.90e{+04}\) \\ \hline
Cranefly &
\(-1.50e{+06}\) &
\(-1.60e{+06}\) &
\(-1.60e{+06}\) &
\(-1.68e{+06}\) \\
& \(-9.28e{+04}\) &
\(-1.01e{+05}\) &
\(-9.60e{+04}\) &
\(-1.04e{+05}\) \\ \hline
Bumblebee &
\(-1.20e{+07}\) &
\(-1.31e{+07}\) &
\(-1.26e{+07}\) &
\(-1.27e{+07}\) \\
& \(-1.02e{+05}\) &
\(-1.12e{+05}\) &
\(-1.07e{+05}\) &
\(-1.08e{+05}\) \\ \hline
Dragonfly &
\(-4.14e{+06}\) &
\(-4.50e{+06}\) &
\(-4.39e{+06}\) &
\(-4.52e{+06}\) \\
& \(-1.41e{+05}\) &
\(-1.54e{+05}\) &
\(-1.60e{+05}\) &
\(-1.82e{+05}\) \\ \hline
Hoverfly &
\(-7.81e{+06}\) &
\(-8.55e{+06}\) &
\(-8.04e{+06}\) &
\(-8.11e{+06}\) \\
& \(-1.20e{+05}\) &
\(-1.31e{+05}\) &
\(-1.24e{+05}\) &
\(-1.25e{+05}\) \\ \hline
Hummingbird &
\(-3.33e{+06}\) &
\(-3.65e{+06}\) &
\(-4.01e{+06}\) &
\(-4.38e{+06}\) \\
& \(-2.00e{+05}\) &
\(-2.14e{+05}\) &
\(-2.40e{+05}\) &
\(-2.62e{+05}\) \\ \hline
\end{tabular}}
\end{table*}

After considering the aforementioned smoothing approximation, we are ready to apply the VOC-based averaging in (\ref{eq:VOC_result}) and analyze stability via linearization and eigenvalues, which in turn characterizes the stability of the natural hovering ES system \eqref{eq:Full_model_ESC}. The reader is referred to the code in \cite{githubdae} for the computations involved in \eqref{eq:VOC_result}, linearization about hovering equilibrium, and eigenvalues computation. However, here we present the end result, provided in Table \ref{tab:eigenvalues}.

Additionally, we made computations and included also in Table \ref{tab:eigenvalues} the eigenvalues of stable open-loop hovering \cite{hassan2018combinedguidancejournal}; the reader is encouraged to recall relevant comments/discussions about open-loop stability provided in Subsection III.A. \textcolor{black}{To obtain the optimal–PID eigenvalues, the averaged model was first derived exactly as in the open-loop case.  
Classical linear-control theory was then applied to incorporate the PID loop and evaluate the closed-loop eigenvalues, following the procedure in \cite{pid_book}.  
All computational steps are available in the accompanying code repository in \cite{githubdae}. It is important to mention that there is a difference in the nature of the stability between PID/ES on one hand and the open-loop case on the other. In the open-loop configuration, the equilibrium need not correspond to hovering (\(w=0\)).  Stability can occur at a non-zero \(w\), leading to steady ascent or descent, as discussed in Section~III.B; only by selecting the exact trim (numerically demanding) does one obtain open-loop hovering.  
By contrast, both the PID loop and the ES controller always stabilize to \(w\!=\!0\), so any disturbance is rejected and altitude is recovered. It can also be noted that in some instances the optimal PID closed-loop is less effective than the open-loop system—in the hawkmoth case, for example, as shown in Table \ref{tab:eigenvalues} (but this should be put in the context that PID stabilization is hovering stabilization, which is not the case in the open-loop case.} Overall, the eigenvalue results confirm the improvement in hovering stability by the natural hovering ES system using different objective functions compared to the open-loop case and \textcolor{black}{optimal PID case}.
The results confirm that via extremum seeking control (ESC), hovering exhibits eigenvalues with more negative real components (i.e., shifts the eigenvalues further into the left-half complex plane),  indicating stronger stability properties and more robustness.

\section{Conclusion and future work} \label{sec:conclusion}
This paper provided a novel characterization (the first of its kind) of the hovering phenomenon in flapping insects and hummingbirds. By describing hovering in flapping insects and hummingbirds as an extremum seeking system (we call it a natural hovering ES system), we were able to \textcolor{black}{start a line of research that may lead to} solving all the puzzle pieces of hovering flight that existed for decades in previous literature. The proposed natural hovering ES system: (i) is very simple computationally; (ii) relies only on the natural oscillation of the wing combined with sensation-based feedback; (iii) real-time; (iv) stable; and (v) model-free (we used available models in literature to test the hypothesis but ES does not need modeling for implementation). We demonstrated the effectiveness of the proposed natural hovering ES via successful simulation trials demonstrating stabilization (even under delays and persistent noise in feedback measurements) in real-time. \textcolor{black}{We also compared the proposed natural hovering ES to open-loop and PID approaches, demonstrating the clear advantages the proposed framework offers.} Additionally, we provided stability analysis, which confirmed the observations provided in the simulations. \textcolor{black}{In the future, we aim at expanding the proposed research to include pitching dynamics, which is at the center of the stability-instability debate and corresponding claims in literature. We also seek experimental validation of natural hovering ES.} 

\newpage
%\clearpage % or \newpage
\bibliography{apssamp.bib}%

\begin{thebibliography}{66}%
\makeatletter
\providecommand \@ifxundefined [1]{%
 \@ifx{#1\undefined}
}%
\providecommand \@ifnum [1]{%
 \ifnum #1\expandafter \@firstoftwo
 \else \expandafter \@secondoftwo
 \fi
}%
\providecommand \@ifx [1]{%
 \ifx #1\expandafter \@firstoftwo
 \else \expandafter \@secondoftwo
 \fi
}%
\providecommand \natexlab [1]{#1}%
\providecommand \enquote  [1]{``#1''}%
\providecommand \bibnamefont  [1]{#1}%
\providecommand \bibfnamefont [1]{#1}%
\providecommand \citenamefont [1]{#1}%
\providecommand \href@noop [0]{\@secondoftwo}%
\providecommand \href [0]{\begingroup \@sanitize@url \@href}%
\providecommand \@href[1]{\@@startlink{#1}\@@href}%
\providecommand \@@href[1]{\endgroup#1\@@endlink}%
\providecommand \@sanitize@url [0]{\catcode `\\12\catcode `\$12\catcode `\&12\catcode `\#12\catcode `\^12\catcode `\_12\catcode `\%12\relax}%
\providecommand \@@startlink[1]{}%
\providecommand \@@endlink[0]{}%
\providecommand \url  [0]{\begingroup\@sanitize@url \@url }%
\providecommand \@url [1]{\endgroup\@href {#1}{\urlprefix }}%
\providecommand \urlprefix  [0]{URL }%
\providecommand \Eprint [0]{\href }%
\providecommand \doibase [0]{https://doi.org/}%
\providecommand \selectlanguage [0]{\@gobble}%
\providecommand \bibinfo  [0]{\@secondoftwo}%
\providecommand \bibfield  [0]{\@secondoftwo}%
\providecommand \translation [1]{[#1]}%
\providecommand \BibitemOpen [0]{}%
\providecommand \bibitemStop [0]{}%
\providecommand \bibitemNoStop [0]{.\EOS\space}%
\providecommand \EOS [0]{\spacefactor3000\relax}%
\providecommand \BibitemShut  [1]{\csname bibitem#1\endcsname}%
\let\auto@bib@innerbib\@empty
%</preamble>
\bibitem [{\citenamefont {Taha}\ \emph {et~al.}(2012)\citenamefont {Taha}, \citenamefont {Hajj},\ and\ \citenamefont {Nayfeh}}]{taha_review}%
  \BibitemOpen
  \bibfield  {author} {\bibinfo {author} {\bibfnamefont {H.~E.}\ \bibnamefont {Taha}}, \bibinfo {author} {\bibfnamefont {M.~R.}\ \bibnamefont {Hajj}},\ and\ \bibinfo {author} {\bibfnamefont {A.~H.}\ \bibnamefont {Nayfeh}},\ }\bibfield  {title} {\bibinfo {title} {Flight dynamics and control of flapping-wing mavs: a review},\ }\href@noop {} {\bibfield  {journal} {\bibinfo  {journal} {Nonlinear Dynamics}\ }\textbf {\bibinfo {volume} {70}},\ \bibinfo {pages} {907} (\bibinfo {year} {2012})}\BibitemShut {NoStop}%
\bibitem [{\citenamefont {Phan}\ and\ \citenamefont {Park}(2019)}]{phan2019insect}%
  \BibitemOpen
  \bibfield  {author} {\bibinfo {author} {\bibfnamefont {H.~V.}\ \bibnamefont {Phan}}\ and\ \bibinfo {author} {\bibfnamefont {H.~C.}\ \bibnamefont {Park}},\ }\bibfield  {title} {\bibinfo {title} {Insect-inspired, tailless, hover-capable flapping-wing robots: Recent progress, challenges, and future directions},\ }\href@noop {} {\bibfield  {journal} {\bibinfo  {journal} {Progress in Aerospace Sciences}\ }\textbf {\bibinfo {volume} {111}},\ \bibinfo {pages} {100573} (\bibinfo {year} {2019})}\BibitemShut {NoStop}%
\bibitem [{\citenamefont {Xuan}\ \emph {et~al.}(2020)\citenamefont {Xuan}, \citenamefont {Hu}, \citenamefont {Yu},\ and\ \citenamefont {Zhang}}]{xuan2020recent}%
  \BibitemOpen
  \bibfield  {author} {\bibinfo {author} {\bibfnamefont {H.}~\bibnamefont {Xuan}}, \bibinfo {author} {\bibfnamefont {J.}~\bibnamefont {Hu}}, \bibinfo {author} {\bibfnamefont {Y.}~\bibnamefont {Yu}},\ and\ \bibinfo {author} {\bibfnamefont {J.}~\bibnamefont {Zhang}},\ }\bibfield  {title} {\bibinfo {title} {Recent progress in aerodynamic modeling methods for flapping flight},\ }\href@noop {} {\bibfield  {journal} {\bibinfo  {journal} {AIP Advances}\ }\textbf {\bibinfo {volume} {10}} (\bibinfo {year} {2020})}\BibitemShut {NoStop}%
\bibitem [{\citenamefont {Weis-Fogh}(1972)}]{weis1972energetics}%
  \BibitemOpen
  \bibfield  {author} {\bibinfo {author} {\bibfnamefont {T.}~\bibnamefont {Weis-Fogh}},\ }\bibfield  {title} {\bibinfo {title} {Energetics of hovering flight in hummingbirds and in drosophila},\ }\href@noop {} {\bibfield  {journal} {\bibinfo  {journal} {Journal of Experimental Biology}\ }\textbf {\bibinfo {volume} {56}},\ \bibinfo {pages} {79} (\bibinfo {year} {1972})}\BibitemShut {NoStop}%
\bibitem [{\citenamefont {Norberg}(1975)}]{norberg1975hovering}%
  \BibitemOpen
  \bibfield  {author} {\bibinfo {author} {\bibfnamefont {R.~{\AA}.}\ \bibnamefont {Norberg}},\ }\bibfield  {title} {\bibinfo {title} {Hovering flight of the dragonfly aeschna juncea l., kinematics and aerodynamics},\ }\href@noop {} {\bibfield  {journal} {\bibinfo  {journal} {Swimming and Flying in Nature: Volume 2}\ ,\ \bibinfo {pages} {763}} (\bibinfo {year} {1975})}\BibitemShut {NoStop}%
\bibitem [{\citenamefont {Dudley}\ and\ \citenamefont {Ellington}(1990{\natexlab{a}})}]{dudley1990mechanics}%
  \BibitemOpen
  \bibfield  {author} {\bibinfo {author} {\bibfnamefont {R.}~\bibnamefont {Dudley}}\ and\ \bibinfo {author} {\bibfnamefont {C.~P.}\ \bibnamefont {Ellington}},\ }\bibfield  {title} {\bibinfo {title} {Mechanics of forward flight in bumblebees: Ii. quasi-steady lift and power requirements},\ }\href@noop {} {\bibfield  {journal} {\bibinfo  {journal} {Journal of Experimental Biology}\ }\textbf {\bibinfo {volume} {148}},\ \bibinfo {pages} {53} (\bibinfo {year} {1990}{\natexlab{a}})}\BibitemShut {NoStop}%
\bibitem [{\citenamefont {Ellington}(1995)}]{ellington1995unsteady}%
  \BibitemOpen
  \bibfield  {author} {\bibinfo {author} {\bibfnamefont {C.}~\bibnamefont {Ellington}},\ }\bibfield  {title} {\bibinfo {title} {Unsteady aerodynamics of insect flight.},\ }in\ \href@noop {} {\emph {\bibinfo {booktitle} {Symposia of the society for experimental biology}}},\ Vol.~\bibinfo {volume} {49}\ (\bibinfo {year} {1995})\ pp.\ \bibinfo {pages} {109--129}\BibitemShut {NoStop}%
\bibitem [{\citenamefont {Willmott}\ and\ \citenamefont {Ellington}(1997)}]{145}%
  \BibitemOpen
  \bibfield  {author} {\bibinfo {author} {\bibfnamefont {A.~P.}\ \bibnamefont {Willmott}}\ and\ \bibinfo {author} {\bibfnamefont {C.~P.}\ \bibnamefont {Ellington}},\ }\bibfield  {title} {\bibinfo {title} {The mechanics of flight in the hawkmoth manduca sexta i. kinematics of hovering and forward flight},\ }\href@noop {} {\bibfield  {journal} {\bibinfo  {journal} {Journal of experimental Biology}\ }\textbf {\bibinfo {volume} {200}},\ \bibinfo {pages} {2705} (\bibinfo {year} {1997})}\BibitemShut {NoStop}%
\bibitem [{\citenamefont {Sun}(2014)}]{sun2014insect}%
  \BibitemOpen
  \bibfield  {author} {\bibinfo {author} {\bibfnamefont {M.}~\bibnamefont {Sun}},\ }\bibfield  {title} {\bibinfo {title} {Insect flight dynamics: stability and control},\ }\href@noop {} {\bibfield  {journal} {\bibinfo  {journal} {Reviews of Modern Physics}\ }\textbf {\bibinfo {volume} {86}},\ \bibinfo {pages} {615} (\bibinfo {year} {2014})}\BibitemShut {NoStop}%
\bibitem [{\citenamefont {Ellington}\ \emph {et~al.}(1996)\citenamefont {Ellington}, \citenamefont {Van Den~Berg}, \citenamefont {Willmott},\ and\ \citenamefont {Thomas}}]{nature1996LEV}%
  \BibitemOpen
  \bibfield  {author} {\bibinfo {author} {\bibfnamefont {C.~P.}\ \bibnamefont {Ellington}}, \bibinfo {author} {\bibfnamefont {C.}~\bibnamefont {Van Den~Berg}}, \bibinfo {author} {\bibfnamefont {A.~P.}\ \bibnamefont {Willmott}},\ and\ \bibinfo {author} {\bibfnamefont {A.~L.}\ \bibnamefont {Thomas}},\ }\bibfield  {title} {\bibinfo {title} {Leading-edge vortices in insect flight},\ }\href@noop {} {\bibfield  {journal} {\bibinfo  {journal} {Nature}\ }\textbf {\bibinfo {volume} {384}},\ \bibinfo {pages} {626} (\bibinfo {year} {1996})}\BibitemShut {NoStop}%
\bibitem [{\citenamefont {Dickinson}\ \emph {et~al.}(1999)\citenamefont {Dickinson}, \citenamefont {Lehmann},\ and\ \citenamefont {Sane}}]{dickinson1999wing}%
  \BibitemOpen
  \bibfield  {author} {\bibinfo {author} {\bibfnamefont {M.~H.}\ \bibnamefont {Dickinson}}, \bibinfo {author} {\bibfnamefont {F.-O.}\ \bibnamefont {Lehmann}},\ and\ \bibinfo {author} {\bibfnamefont {S.~P.}\ \bibnamefont {Sane}},\ }\bibfield  {title} {\bibinfo {title} {Wing rotation and the aerodynamic basis of insect flight},\ }\href@noop {} {\bibfield  {journal} {\bibinfo  {journal} {Science}\ }\textbf {\bibinfo {volume} {284}},\ \bibinfo {pages} {1954} (\bibinfo {year} {1999})}\BibitemShut {NoStop}%
\bibitem [{\citenamefont {Taylor}\ and\ \citenamefont {Thomas}(2003)}]{taylor2003dynamic}%
  \BibitemOpen
  \bibfield  {author} {\bibinfo {author} {\bibfnamefont {G.~K.}\ \bibnamefont {Taylor}}\ and\ \bibinfo {author} {\bibfnamefont {A.~L.}\ \bibnamefont {Thomas}},\ }\bibfield  {title} {\bibinfo {title} {Dynamic flight stability in the desert locust schistocerca gregaria},\ }\href@noop {} {\bibfield  {journal} {\bibinfo  {journal} {Journal of Experimental Biology}\ }\textbf {\bibinfo {volume} {206}},\ \bibinfo {pages} {2803} (\bibinfo {year} {2003})}\BibitemShut {NoStop}%
\bibitem [{\citenamefont {Taylor}\ and\ \citenamefont {{\.Z}bikowski}(2005)}]{taylor2005nonlinear}%
  \BibitemOpen
  \bibfield  {author} {\bibinfo {author} {\bibfnamefont {G.~K.}\ \bibnamefont {Taylor}}\ and\ \bibinfo {author} {\bibfnamefont {R.}~\bibnamefont {{\.Z}bikowski}},\ }\bibfield  {title} {\bibinfo {title} {Nonlinear time-periodic models of the longitudinal flight dynamics of desert locusts schistocerca gregaria},\ }\href@noop {} {\bibfield  {journal} {\bibinfo  {journal} {Journal of the Royal Society Interface}\ }\textbf {\bibinfo {volume} {2}},\ \bibinfo {pages} {197} (\bibinfo {year} {2005})}\BibitemShut {NoStop}%
\bibitem [{\citenamefont {Zhang}\ \emph {et~al.}(2012)\citenamefont {Zhang}, \citenamefont {Wu},\ and\ \citenamefont {Sun}}]{zhang2012lateral}%
  \BibitemOpen
  \bibfield  {author} {\bibinfo {author} {\bibfnamefont {Y.-L.}\ \bibnamefont {Zhang}}, \bibinfo {author} {\bibfnamefont {J.-H.}\ \bibnamefont {Wu}},\ and\ \bibinfo {author} {\bibfnamefont {M.}~\bibnamefont {Sun}},\ }\bibfield  {title} {\bibinfo {title} {Lateral dynamic flight stability of hovering insects: theory vs. numerical simulation},\ }\href@noop {} {\bibfield  {journal} {\bibinfo  {journal} {Acta Mechanica Sinica}\ }\textbf {\bibinfo {volume} {28}},\ \bibinfo {pages} {221} (\bibinfo {year} {2012})}\BibitemShut {NoStop}%
\bibitem [{\citenamefont {Xu}\ and\ \citenamefont {Sun}(2013)}]{xu2013lateral}%
  \BibitemOpen
  \bibfield  {author} {\bibinfo {author} {\bibfnamefont {N.}~\bibnamefont {Xu}}\ and\ \bibinfo {author} {\bibfnamefont {M.}~\bibnamefont {Sun}},\ }\bibfield  {title} {\bibinfo {title} {Lateral dynamic flight stability of a model bumblebee in hovering and forward flight},\ }\href@noop {} {\bibfield  {journal} {\bibinfo  {journal} {Journal of theoretical biology}\ }\textbf {\bibinfo {volume} {319}},\ \bibinfo {pages} {102} (\bibinfo {year} {2013})}\BibitemShut {NoStop}%
\bibitem [{\citenamefont {Kar{\'a}sek}(2014)}]{DelftDissertation2014}%
  \BibitemOpen
  \bibfield  {author} {\bibinfo {author} {\bibfnamefont {M.}~\bibnamefont {Kar{\'a}sek}},\ }\bibfield  {title} {\bibinfo {title} {Robotic hummingbird: Design of a control mechanism for a hovering flapping wing micro air vehicle},\ }\href@noop {} {\bibfield  {journal} {\bibinfo  {journal} {Universite libre de Bruxelles: Bruxelles, Belgium}\ } (\bibinfo {year} {2014})}\BibitemShut {NoStop}%
\bibitem [{\citenamefont {Ristroph}\ \emph {et~al.}(2010)\citenamefont {Ristroph}, \citenamefont {Bergou}, \citenamefont {Ristroph}, \citenamefont {Coumes}, \citenamefont {Berman}, \citenamefont {Guckenheimer}, \citenamefont {Wang},\ and\ \citenamefont {Cohen}}]{feedback25ristroph2010discovering}%
  \BibitemOpen
  \bibfield  {author} {\bibinfo {author} {\bibfnamefont {L.}~\bibnamefont {Ristroph}}, \bibinfo {author} {\bibfnamefont {A.~J.}\ \bibnamefont {Bergou}}, \bibinfo {author} {\bibfnamefont {G.}~\bibnamefont {Ristroph}}, \bibinfo {author} {\bibfnamefont {K.}~\bibnamefont {Coumes}}, \bibinfo {author} {\bibfnamefont {G.~J.}\ \bibnamefont {Berman}}, \bibinfo {author} {\bibfnamefont {J.}~\bibnamefont {Guckenheimer}}, \bibinfo {author} {\bibfnamefont {Z.~J.}\ \bibnamefont {Wang}},\ and\ \bibinfo {author} {\bibfnamefont {I.}~\bibnamefont {Cohen}},\ }\bibfield  {title} {\bibinfo {title} {Discovering the flight autostabilizer of fruit flies by inducing aerial stumbles},\ }\href@noop {} {\bibfield  {journal} {\bibinfo  {journal} {Proceedings of the National Academy of Sciences}\ }\textbf {\bibinfo {volume} {107}},\ \bibinfo {pages} {4820} (\bibinfo {year} {2010})}\BibitemShut {NoStop}%
\bibitem [{\citenamefont {Lyu}\ and\ \citenamefont {Sun}(2022)}]{lyu2022dynamic}%
  \BibitemOpen
  \bibfield  {author} {\bibinfo {author} {\bibfnamefont {Y.~Z.}\ \bibnamefont {Lyu}}\ and\ \bibinfo {author} {\bibfnamefont {M.}~\bibnamefont {Sun}},\ }\bibfield  {title} {\bibinfo {title} {Dynamic stability in hovering flight of insects with different sizes},\ }\href@noop {} {\bibfield  {journal} {\bibinfo  {journal} {Physical Review E}\ }\textbf {\bibinfo {volume} {105}},\ \bibinfo {pages} {054403} (\bibinfo {year} {2022})}\BibitemShut {NoStop}%
\bibitem [{\citenamefont {Taha}\ \emph {et~al.}(2015)\citenamefont {Taha}, \citenamefont {Tahmasian}, \citenamefont {Woolsey}, \citenamefont {Nayfeh},\ and\ \citenamefont {Hajj}}]{taha2015need}%
  \BibitemOpen
  \bibfield  {author} {\bibinfo {author} {\bibfnamefont {H.~E.}\ \bibnamefont {Taha}}, \bibinfo {author} {\bibfnamefont {S.}~\bibnamefont {Tahmasian}}, \bibinfo {author} {\bibfnamefont {C.~A.}\ \bibnamefont {Woolsey}}, \bibinfo {author} {\bibfnamefont {A.~H.}\ \bibnamefont {Nayfeh}},\ and\ \bibinfo {author} {\bibfnamefont {M.~R.}\ \bibnamefont {Hajj}},\ }\bibfield  {title} {\bibinfo {title} {The need for higher-order averaging in the stability analysis of hovering, flapping-wing flight},\ }\href@noop {} {\bibfield  {journal} {\bibinfo  {journal} {Bioinspiration \& biomimetics}\ }\textbf {\bibinfo {volume} {10}},\ \bibinfo {pages} {016002} (\bibinfo {year} {2015})}\BibitemShut {NoStop}%
\bibitem [{\citenamefont {Taha}\ \emph {et~al.}(2020)\citenamefont {Taha}, \citenamefont {Kiani}, \citenamefont {Hedrick},\ and\ \citenamefont {Greeter}}]{taha2020vibrational}%
  \BibitemOpen
  \bibfield  {author} {\bibinfo {author} {\bibfnamefont {H.~E.}\ \bibnamefont {Taha}}, \bibinfo {author} {\bibfnamefont {M.}~\bibnamefont {Kiani}}, \bibinfo {author} {\bibfnamefont {T.~L.}\ \bibnamefont {Hedrick}},\ and\ \bibinfo {author} {\bibfnamefont {J.~S.}\ \bibnamefont {Greeter}},\ }\bibfield  {title} {\bibinfo {title} {Vibrational control: A hidden stabilization mechanism in insect flight},\ }\href@noop {} {\bibfield  {journal} {\bibinfo  {journal} {Science robotics}\ }\textbf {\bibinfo {volume} {5}},\ \bibinfo {pages} {eabb1502} (\bibinfo {year} {2020})}\BibitemShut {NoStop}%
\bibitem [{\citenamefont {Fuller}\ \emph {et~al.}(2014)\citenamefont {Fuller}, \citenamefont {Straw}, \citenamefont {Peek}, \citenamefont {Murray},\ and\ \citenamefont {Dickinson}}]{fuller2014flying}%
  \BibitemOpen
  \bibfield  {author} {\bibinfo {author} {\bibfnamefont {S.~B.}\ \bibnamefont {Fuller}}, \bibinfo {author} {\bibfnamefont {A.~D.}\ \bibnamefont {Straw}}, \bibinfo {author} {\bibfnamefont {M.~Y.}\ \bibnamefont {Peek}}, \bibinfo {author} {\bibfnamefont {R.~M.}\ \bibnamefont {Murray}},\ and\ \bibinfo {author} {\bibfnamefont {M.~H.}\ \bibnamefont {Dickinson}},\ }\bibfield  {title} {\bibinfo {title} {Flying drosophila stabilize their vision-based velocity controller by sensing wind with their antennae},\ }\href@noop {} {\bibfield  {journal} {\bibinfo  {journal} {Proceedings of the National Academy of Sciences}\ }\textbf {\bibinfo {volume} {111}},\ \bibinfo {pages} {E1182} (\bibinfo {year} {2014})}\BibitemShut {NoStop}%
\bibitem [{\citenamefont {Ristroph}\ \emph {et~al.}(2013)\citenamefont {Ristroph}, \citenamefont {Ristroph}, \citenamefont {Morozova}, \citenamefont {Bergou}, \citenamefont {Chang}, \citenamefont {Guckenheimer}, \citenamefont {Wang},\ and\ \citenamefont {Cohen}}]{sensorroyalristroph2013active}%
  \BibitemOpen
  \bibfield  {author} {\bibinfo {author} {\bibfnamefont {L.}~\bibnamefont {Ristroph}}, \bibinfo {author} {\bibfnamefont {G.}~\bibnamefont {Ristroph}}, \bibinfo {author} {\bibfnamefont {S.}~\bibnamefont {Morozova}}, \bibinfo {author} {\bibfnamefont {A.~J.}\ \bibnamefont {Bergou}}, \bibinfo {author} {\bibfnamefont {S.}~\bibnamefont {Chang}}, \bibinfo {author} {\bibfnamefont {J.}~\bibnamefont {Guckenheimer}}, \bibinfo {author} {\bibfnamefont {Z.~J.}\ \bibnamefont {Wang}},\ and\ \bibinfo {author} {\bibfnamefont {I.}~\bibnamefont {Cohen}},\ }\bibfield  {title} {\bibinfo {title} {Active and passive stabilization of body pitch in insect flight},\ }\href@noop {} {\bibfield  {journal} {\bibinfo  {journal} {Journal of The Royal Society Interface}\ }\textbf {\bibinfo {volume} {10}},\ \bibinfo {pages} {20130237} (\bibinfo {year} {2013})}\BibitemShut {NoStop}%
\bibitem [{\citenamefont {Cheng}\ \emph {et~al.}(2011)\citenamefont {Cheng}, \citenamefont {Deng},\ and\ \citenamefont {Hedrick}}]{feedback26cheng2011mechanics}%
  \BibitemOpen
  \bibfield  {author} {\bibinfo {author} {\bibfnamefont {B.}~\bibnamefont {Cheng}}, \bibinfo {author} {\bibfnamefont {X.}~\bibnamefont {Deng}},\ and\ \bibinfo {author} {\bibfnamefont {T.~L.}\ \bibnamefont {Hedrick}},\ }\bibfield  {title} {\bibinfo {title} {The mechanics and control of pitching manoeuvres in a freely flying hawkmoth (manduca sexta)},\ }\href@noop {} {\bibfield  {journal} {\bibinfo  {journal} {Journal of Experimental Biology}\ }\textbf {\bibinfo {volume} {214}},\ \bibinfo {pages} {4092} (\bibinfo {year} {2011})}\BibitemShut {NoStop}%
\bibitem [{\citenamefont {Taylor}\ and\ \citenamefont {Krapp}(2007)}]{sensors1taylor2007sensory}%
  \BibitemOpen
  \bibfield  {author} {\bibinfo {author} {\bibfnamefont {G.~K.}\ \bibnamefont {Taylor}}\ and\ \bibinfo {author} {\bibfnamefont {H.~G.}\ \bibnamefont {Krapp}},\ }\bibfield  {title} {\bibinfo {title} {Sensory systems and flight stability: what do insects measure and why?},\ }\href@noop {} {\bibfield  {journal} {\bibinfo  {journal} {Advances in insect physiology}\ }\textbf {\bibinfo {volume} {34}},\ \bibinfo {pages} {231} (\bibinfo {year} {2007})}\BibitemShut {NoStop}%
\bibitem [{\citenamefont {Deng}\ \emph {et~al.}(2006{\natexlab{a}})\citenamefont {Deng}, \citenamefont {Schenato}, \citenamefont {Wu},\ and\ \citenamefont {Sastry}}]{IEEETransaction}%
  \BibitemOpen
  \bibfield  {author} {\bibinfo {author} {\bibfnamefont {X.}~\bibnamefont {Deng}}, \bibinfo {author} {\bibfnamefont {L.}~\bibnamefont {Schenato}}, \bibinfo {author} {\bibfnamefont {W.~C.}\ \bibnamefont {Wu}},\ and\ \bibinfo {author} {\bibfnamefont {S.}~\bibnamefont {Sastry}},\ }\bibfield  {title} {\bibinfo {title} {Flapping flight for biomimetic robotic insects: part i-system modeling},\ }\href {https://doi.org/10.1109/TRO.2006.875480} {\bibfield  {journal} {\bibinfo  {journal} {IEEE Transactions on Robotics}\ }\textbf {\bibinfo {volume} {22}},\ \bibinfo {pages} {776} (\bibinfo {year} {2006}{\natexlab{a}})}\BibitemShut {NoStop}%
\bibitem [{\citenamefont {Rapp}\ and\ \citenamefont {Nawrot}(2020)}]{rapp2020spiking}%
  \BibitemOpen
  \bibfield  {author} {\bibinfo {author} {\bibfnamefont {H.}~\bibnamefont {Rapp}}\ and\ \bibinfo {author} {\bibfnamefont {M.~P.}\ \bibnamefont {Nawrot}},\ }\bibfield  {title} {\bibinfo {title} {A spiking neural program for sensorimotor control during foraging in flying insects},\ }\href@noop {} {\bibfield  {journal} {\bibinfo  {journal} {Proceedings of the National Academy of Sciences}\ }\textbf {\bibinfo {volume} {117}},\ \bibinfo {pages} {28412} (\bibinfo {year} {2020})}\BibitemShut {NoStop}%
\bibitem [{\citenamefont {Deng}\ \emph {et~al.}(2006{\natexlab{b}})\citenamefont {Deng}, \citenamefont {Schenato},\ and\ \citenamefont {Sastry}}]{24}%
  \BibitemOpen
  \bibfield  {author} {\bibinfo {author} {\bibfnamefont {X.}~\bibnamefont {Deng}}, \bibinfo {author} {\bibfnamefont {L.}~\bibnamefont {Schenato}},\ and\ \bibinfo {author} {\bibfnamefont {S.~S.}\ \bibnamefont {Sastry}},\ }\bibfield  {title} {\bibinfo {title} {Flapping flight for biomimetic robotic insects: Part ii-flight control design},\ }\href@noop {} {\bibfield  {journal} {\bibinfo  {journal} {IEEE transactions on robotics}\ }\textbf {\bibinfo {volume} {22}},\ \bibinfo {pages} {789} (\bibinfo {year} {2006}{\natexlab{b}})}\BibitemShut {NoStop}%
\bibitem [{\citenamefont {Kar{\'a}sek}\ and\ \citenamefont {Preumont}(2012)}]{Delft_Conference}%
  \BibitemOpen
  \bibfield  {author} {\bibinfo {author} {\bibfnamefont {M.}~\bibnamefont {Kar{\'a}sek}}\ and\ \bibinfo {author} {\bibfnamefont {A.}~\bibnamefont {Preumont}},\ }\bibfield  {title} {\bibinfo {title} {Simulation of flight control of a hummingbird like robot near hover},\ }\href@noop {} {\bibfield  {journal} {\bibinfo  {journal} {Engineering Mechanics}\ }\textbf {\bibinfo {volume} {322}} (\bibinfo {year} {2012})}\BibitemShut {NoStop}%
\bibitem [{\citenamefont {Serrani}\ \emph {et~al.}(2010)\citenamefont {Serrani}, \citenamefont {Keller}, \citenamefont {Bolender},\ and\ \citenamefont {Doman}}]{serrani2010robustAIAA}%
  \BibitemOpen
  \bibfield  {author} {\bibinfo {author} {\bibfnamefont {A.}~\bibnamefont {Serrani}}, \bibinfo {author} {\bibfnamefont {B.}~\bibnamefont {Keller}}, \bibinfo {author} {\bibfnamefont {M.}~\bibnamefont {Bolender}},\ and\ \bibinfo {author} {\bibfnamefont {D.}~\bibnamefont {Doman}},\ }\bibfield  {title} {\bibinfo {title} {Robust control of a 3-dof flapping wing micro air vehicle},\ }in\ \href@noop {} {\emph {\bibinfo {booktitle} {AIAA Guidance, Navigation, and Control Conference}}}\ (\bibinfo {year} {2010})\ p.\ \bibinfo {pages} {7709}\BibitemShut {NoStop}%
\bibitem [{\citenamefont {Doman}\ \emph {et~al.}(2009)\citenamefont {Doman}, \citenamefont {Oppenheimer},\ and\ \citenamefont {Sigthorsson}}]{doman2009dynamics}%
  \BibitemOpen
  \bibfield  {author} {\bibinfo {author} {\bibfnamefont {D.}~\bibnamefont {Doman}}, \bibinfo {author} {\bibfnamefont {M.}~\bibnamefont {Oppenheimer}},\ and\ \bibinfo {author} {\bibfnamefont {D.}~\bibnamefont {Sigthorsson}},\ }\bibfield  {title} {\bibinfo {title} {Dynamics and control of a minimally actuated biomimetic vehicle: Part i-aerodynamic model},\ }in\ \href@noop {} {\emph {\bibinfo {booktitle} {AIAA guidance, navigation, and control conference}}}\ (\bibinfo {year} {2009})\ p.\ \bibinfo {pages} {6160}\BibitemShut {NoStop}%
\bibitem [{\citenamefont {Oppenheimer}\ \emph {et~al.}(2009)\citenamefont {Oppenheimer}, \citenamefont {Doman},\ and\ \citenamefont {Sigthorsson}}]{oppenheimer2009dynamics}%
  \BibitemOpen
  \bibfield  {author} {\bibinfo {author} {\bibfnamefont {M.}~\bibnamefont {Oppenheimer}}, \bibinfo {author} {\bibfnamefont {D.}~\bibnamefont {Doman}},\ and\ \bibinfo {author} {\bibfnamefont {D.}~\bibnamefont {Sigthorsson}},\ }\bibfield  {title} {\bibinfo {title} {Dynamics and control of a minimally actuated biomimetic vehicle: Part ii-control},\ }in\ \href@noop {} {\emph {\bibinfo {booktitle} {AIAA Guidance, Navigation, and Control Conference}}}\ (\bibinfo {year} {2009})\ p.\ \bibinfo {pages} {6161}\BibitemShut {NoStop}%
\bibitem [{\citenamefont {Sun}\ and\ \citenamefont {Lan}(2004)}]{sun2004computational}%
  \BibitemOpen
  \bibfield  {author} {\bibinfo {author} {\bibfnamefont {M.}~\bibnamefont {Sun}}\ and\ \bibinfo {author} {\bibfnamefont {S.~L.}\ \bibnamefont {Lan}},\ }\bibfield  {title} {\bibinfo {title} {A computational study of the aerodynamic forces and power requirements of dragonfly (aeschna juncea) hovering},\ }\href@noop {} {\bibfield  {journal} {\bibinfo  {journal} {Journal of Experimental Biology}\ }\textbf {\bibinfo {volume} {207}},\ \bibinfo {pages} {1887} (\bibinfo {year} {2004})}\BibitemShut {NoStop}%
\bibitem [{\citenamefont {Hedrick}\ \emph {et~al.}(2015)\citenamefont {Hedrick}, \citenamefont {Combes},\ and\ \citenamefont {Miller}}]{hedrick2015recent}%
  \BibitemOpen
  \bibfield  {author} {\bibinfo {author} {\bibfnamefont {T.~L.}\ \bibnamefont {Hedrick}}, \bibinfo {author} {\bibfnamefont {S.~A.}\ \bibnamefont {Combes}},\ and\ \bibinfo {author} {\bibfnamefont {L.~A.}\ \bibnamefont {Miller}},\ }\bibfield  {title} {\bibinfo {title} {Recent developments in the study of insect flight},\ }\href@noop {} {\bibfield  {journal} {\bibinfo  {journal} {Canadian Journal of Zoology}\ }\textbf {\bibinfo {volume} {93}},\ \bibinfo {pages} {925} (\bibinfo {year} {2015})}\BibitemShut {NoStop}%
\bibitem [{\citenamefont {Ariyur}\ and\ \citenamefont {Krstic}(2003)}]{ariyur2003real}%
  \BibitemOpen
  \bibfield  {author} {\bibinfo {author} {\bibfnamefont {K.~B.}\ \bibnamefont {Ariyur}}\ and\ \bibinfo {author} {\bibfnamefont {M.}~\bibnamefont {Krstic}},\ }\href@noop {} {\emph {\bibinfo {title} {Real-time optimization by extremum-seeking control}}}\ (\bibinfo  {publisher} {John Wiley \& Sons},\ \bibinfo {year} {2003})\BibitemShut {NoStop}%
\bibitem [{\citenamefont {Scheinker}\ and\ \citenamefont {Krsti{\'c}}(2017)}]{scheinker2017model}%
  \BibitemOpen
  \bibfield  {author} {\bibinfo {author} {\bibfnamefont {A.}~\bibnamefont {Scheinker}}\ and\ \bibinfo {author} {\bibfnamefont {M.}~\bibnamefont {Krsti{\'c}}},\ }\href@noop {} {\emph {\bibinfo {title} {Model-free stabilization by extremum seeking}}}\ (\bibinfo  {publisher} {Springer},\ \bibinfo {year} {2017})\BibitemShut {NoStop}%
\bibitem [{\citenamefont {Scheinker}(2024)}]{scheinker2024100}%
  \BibitemOpen
  \bibfield  {author} {\bibinfo {author} {\bibfnamefont {A.}~\bibnamefont {Scheinker}},\ }\bibfield  {title} {\bibinfo {title} {100 years of extremum seeking: A survey},\ }\href@noop {} {\bibfield  {journal} {\bibinfo  {journal} {Automatica}\ }\textbf {\bibinfo {volume} {161}},\ \bibinfo {pages} {111481} (\bibinfo {year} {2024})}\BibitemShut {NoStop}%
\bibitem [{\citenamefont {Cochran}\ \emph {et~al.}(2009)\citenamefont {Cochran}, \citenamefont {Kanso}, \citenamefont {Kelly}, \citenamefont {Xiong},\ and\ \citenamefont {Krstic}}]{ESCfishcochran2009source}%
  \BibitemOpen
  \bibfield  {author} {\bibinfo {author} {\bibfnamefont {J.}~\bibnamefont {Cochran}}, \bibinfo {author} {\bibfnamefont {E.}~\bibnamefont {Kanso}}, \bibinfo {author} {\bibfnamefont {S.~D.}\ \bibnamefont {Kelly}}, \bibinfo {author} {\bibfnamefont {H.}~\bibnamefont {Xiong}},\ and\ \bibinfo {author} {\bibfnamefont {M.}~\bibnamefont {Krstic}},\ }\bibfield  {title} {\bibinfo {title} {Source seeking for two nonholonomic models of fish locomotion},\ }\href@noop {} {\bibfield  {journal} {\bibinfo  {journal} {IEEE Transactions on Robotics}\ }\textbf {\bibinfo {volume} {25}},\ \bibinfo {pages} {1166} (\bibinfo {year} {2009})}\BibitemShut {NoStop}%
\bibitem [{\citenamefont {Krstic}\ and\ \citenamefont {Cochran}(2008)}]{krstic2008extremum}%
  \BibitemOpen
  \bibfield  {author} {\bibinfo {author} {\bibfnamefont {M.}~\bibnamefont {Krstic}}\ and\ \bibinfo {author} {\bibfnamefont {J.}~\bibnamefont {Cochran}},\ }\bibfield  {title} {\bibinfo {title} {Extremum seeking for motion optimization: From bacteria to nonholonomic vehicles},\ }in\ \href@noop {} {\emph {\bibinfo {booktitle} {2008 Chinese Control and Decision Conference}}}\ (\bibinfo {organization} {IEEE},\ \bibinfo {year} {2008})\ pp.\ \bibinfo {pages} {18--27}\BibitemShut {NoStop}%
\bibitem [{\citenamefont {Abdelgalil}\ \emph {et~al.}(2022)\citenamefont {Abdelgalil}, \citenamefont {Aboelkassem},\ and\ \citenamefont {Taha}}]{abdelgalil2022sea}%
  \BibitemOpen
  \bibfield  {author} {\bibinfo {author} {\bibfnamefont {M.}~\bibnamefont {Abdelgalil}}, \bibinfo {author} {\bibfnamefont {Y.}~\bibnamefont {Aboelkassem}},\ and\ \bibinfo {author} {\bibfnamefont {H.}~\bibnamefont {Taha}},\ }\bibfield  {title} {\bibinfo {title} {Sea urchin sperm exploit extremum seeking control to find the egg},\ }\href@noop {} {\bibfield  {journal} {\bibinfo  {journal} {Physical Review E}\ }\textbf {\bibinfo {volume} {106}},\ \bibinfo {pages} {L062401} (\bibinfo {year} {2022})}\BibitemShut {NoStop}%
\bibitem [{\citenamefont {Pokhrel}\ and\ \citenamefont {Eisa}(2022)}]{pokhrel2022novel}%
  \BibitemOpen
  \bibfield  {author} {\bibinfo {author} {\bibfnamefont {S.}~\bibnamefont {Pokhrel}}\ and\ \bibinfo {author} {\bibfnamefont {S.~A.}\ \bibnamefont {Eisa}},\ }\bibfield  {title} {\bibinfo {title} {A novel hypothesis for how albatrosses optimize their flight physics in real-time: an extremum seeking model and control for dynamic soaring},\ }\href@noop {} {\bibfield  {journal} {\bibinfo  {journal} {Bioinspiration \& Biomimetics}\ }\textbf {\bibinfo {volume} {18}},\ \bibinfo {pages} {016014} (\bibinfo {year} {2022})}\BibitemShut {NoStop}%
\bibitem [{\citenamefont {Eisa}\ and\ \citenamefont {Pokhrel}(2023)}]{eisa2023analyzing}%
  \BibitemOpen
  \bibfield  {author} {\bibinfo {author} {\bibfnamefont {S.~A.}\ \bibnamefont {Eisa}}\ and\ \bibinfo {author} {\bibfnamefont {S.}~\bibnamefont {Pokhrel}},\ }\bibfield  {title} {\bibinfo {title} {Analyzing and mimicking the optimized flight physics of soaring birds: A differential geometric control and extremum seeking system approach with real time implementation},\ }\href@noop {} {\bibfield  {journal} {\bibinfo  {journal} {SIAM Journal on Applied Mathematics}\ ,\ \bibinfo {pages} {S82}} (\bibinfo {year} {2023})}\BibitemShut {NoStop}%
\bibitem [{\citenamefont {Bullo}(2002)}]{bullo2002averaging}%
  \BibitemOpen
  \bibfield  {author} {\bibinfo {author} {\bibfnamefont {F.}~\bibnamefont {Bullo}},\ }\bibfield  {title} {\bibinfo {title} {Averaging and vibrational control of mechanical systems},\ }\href@noop {} {\bibfield  {journal} {\bibinfo  {journal} {SIAM Journal on Control and Optimization}\ }\textbf {\bibinfo {volume} {41}},\ \bibinfo {pages} {542} (\bibinfo {year} {2002})}\BibitemShut {NoStop}%
\bibitem [{\citenamefont {Pokhrel}\ and\ \citenamefont {Eisa}(2023)}]{pokhrel2023higher}%
  \BibitemOpen
  \bibfield  {author} {\bibinfo {author} {\bibfnamefont {S.}~\bibnamefont {Pokhrel}}\ and\ \bibinfo {author} {\bibfnamefont {S.~A.}\ \bibnamefont {Eisa}},\ }\bibfield  {title} {\bibinfo {title} {Higher-order lie bracket approximation and averaging of control-affine systems with application to extremum seeking},\ }\href@noop {} {\bibfield  {journal} {\bibinfo  {journal} {arXiv preprint arXiv:2310.07092}\ } (\bibinfo {year} {2023})}\BibitemShut {NoStop}%
\bibitem [{\citenamefont {Elgohary}\ and\ \citenamefont {Eisa}(2025{\natexlab{a}})}]{elgohary2025extremum}%
  \BibitemOpen
  \bibfield  {author} {\bibinfo {author} {\bibfnamefont {A.~A.}\ \bibnamefont {Elgohary}}\ and\ \bibinfo {author} {\bibfnamefont {S.~A.}\ \bibnamefont {Eisa}},\ }\bibfield  {title} {\bibinfo {title} {Extremum seeking for controlled vibrational stabilization of mechanical systems: A variation-of-constant averaging approach inspired by flapping insects mechanics},\ }\href@noop {} {\bibfield  {journal} {\bibinfo  {journal} {IEEE Control Systems Letters}\ } (\bibinfo {year} {2025}{\natexlab{a}})}\BibitemShut {NoStop}%
\bibitem [{\citenamefont {Taha}\ \emph {et~al.}(2014)\citenamefont {Taha}, \citenamefont {Hajj},\ and\ \citenamefont {Nayfeh}}]{taha2014longitudinalmodelguidance}%
  \BibitemOpen
  \bibfield  {author} {\bibinfo {author} {\bibfnamefont {H.~E.}\ \bibnamefont {Taha}}, \bibinfo {author} {\bibfnamefont {M.~R.}\ \bibnamefont {Hajj}},\ and\ \bibinfo {author} {\bibfnamefont {A.~H.}\ \bibnamefont {Nayfeh}},\ }\bibfield  {title} {\bibinfo {title} {Longitudinal flight dynamics of hovering mavs/insects},\ }\href@noop {} {\bibfield  {journal} {\bibinfo  {journal} {Journal of Guidance, Control, and Dynamics}\ }\textbf {\bibinfo {volume} {37}},\ \bibinfo {pages} {970} (\bibinfo {year} {2014})}\BibitemShut {NoStop}%
\bibitem [{\citenamefont {Hassan}\ and\ \citenamefont {Taha}(2018)}]{hassan2018combinedguidancejournal}%
  \BibitemOpen
  \bibfield  {author} {\bibinfo {author} {\bibfnamefont {A.~M.}\ \bibnamefont {Hassan}}\ and\ \bibinfo {author} {\bibfnamefont {H.~E.}\ \bibnamefont {Taha}},\ }\bibfield  {title} {\bibinfo {title} {Combined averaging--shooting approach for the analysis of flapping flight dynamics},\ }\href@noop {} {\bibfield  {journal} {\bibinfo  {journal} {Journal of guidance, control, and dynamics}\ }\textbf {\bibinfo {volume} {41}},\ \bibinfo {pages} {542} (\bibinfo {year} {2018})}\BibitemShut {NoStop}%
\bibitem [{\citenamefont {Taha}\ and\ \citenamefont {Hajj}(2013)}]{taha2013unsteady}%
  \BibitemOpen
  \bibfield  {author} {\bibinfo {author} {\bibfnamefont {H.}~\bibnamefont {Taha}}\ and\ \bibinfo {author} {\bibfnamefont {M.}~\bibnamefont {Hajj}},\ }\bibfield  {title} {\bibinfo {title} {Unsteady nonlinear aerodynamics of hovering mavs/insects},\ }in\ \href@noop {} {\emph {\bibinfo {booktitle} {51st AIAA Aerospace Sciences Meeting Including the New Horizons Forum and Aerospace Exposition}}}\ (\bibinfo {year} {2013})\ p.\ \bibinfo {pages} {504}\BibitemShut {NoStop}%
\bibitem [{\citenamefont {Wang}\ \emph {et~al.}(2004)\citenamefont {Wang}, \citenamefont {Birch},\ and\ \citenamefont {Dickinson}}]{141wang2004unsteady}%
  \BibitemOpen
  \bibfield  {author} {\bibinfo {author} {\bibfnamefont {Z.~J.}\ \bibnamefont {Wang}}, \bibinfo {author} {\bibfnamefont {J.~M.}\ \bibnamefont {Birch}},\ and\ \bibinfo {author} {\bibfnamefont {M.~H.}\ \bibnamefont {Dickinson}},\ }\bibfield  {title} {\bibinfo {title} {Unsteady forces and flows in low reynolds number hovering flight: two-dimensional computations vs robotic wing experiments},\ }\href@noop {} {\bibfield  {journal} {\bibinfo  {journal} {Journal of Experimental Biology}\ }\textbf {\bibinfo {volume} {207}},\ \bibinfo {pages} {449} (\bibinfo {year} {2004})}\BibitemShut {NoStop}%
\bibitem [{\citenamefont {Andersen}\ \emph {et~al.}(2005)\citenamefont {Andersen}, \citenamefont {Pesavento},\ and\ \citenamefont {Wang}}]{1andersen2005unsteady}%
  \BibitemOpen
  \bibfield  {author} {\bibinfo {author} {\bibfnamefont {A.}~\bibnamefont {Andersen}}, \bibinfo {author} {\bibfnamefont {U.}~\bibnamefont {Pesavento}},\ and\ \bibinfo {author} {\bibfnamefont {Z.~J.}\ \bibnamefont {Wang}},\ }\bibfield  {title} {\bibinfo {title} {Unsteady aerodynamics of fluttering and tumbling plates},\ }\href@noop {} {\bibfield  {journal} {\bibinfo  {journal} {Journal of Fluid Mechanics}\ }\textbf {\bibinfo {volume} {541}},\ \bibinfo {pages} {65} (\bibinfo {year} {2005})}\BibitemShut {NoStop}%
\bibitem [{\citenamefont {Taha}\ \emph {et~al.}(2016)\citenamefont {Taha}, \citenamefont {Woolsey},\ and\ \citenamefont {Hajj}}]{taha2016geometriclong}%
  \BibitemOpen
  \bibfield  {author} {\bibinfo {author} {\bibfnamefont {H.~E.}\ \bibnamefont {Taha}}, \bibinfo {author} {\bibfnamefont {C.~A.}\ \bibnamefont {Woolsey}},\ and\ \bibinfo {author} {\bibfnamefont {M.~R.}\ \bibnamefont {Hajj}},\ }\bibfield  {title} {\bibinfo {title} {Geometric control approach to longitudinal stability of flapping flight},\ }\href@noop {} {\bibfield  {journal} {\bibinfo  {journal} {Journal of Guidance, Control, and Dynamics}\ }\textbf {\bibinfo {volume} {39}},\ \bibinfo {pages} {214} (\bibinfo {year} {2016})}\BibitemShut {NoStop}%
\bibitem [{\citenamefont {Hassan}\ and\ \citenamefont {Taha}(2017)}]{hassan2017combinedscitech}%
  \BibitemOpen
  \bibfield  {author} {\bibinfo {author} {\bibfnamefont {A.~M.}\ \bibnamefont {Hassan}}\ and\ \bibinfo {author} {\bibfnamefont {H.~E.}\ \bibnamefont {Taha}},\ }\bibfield  {title} {\bibinfo {title} {A combined averaging-shooting approach for the trim analysis of hovering insects/flapping-wing micro-air-vehicles},\ }in\ \href@noop {} {\emph {\bibinfo {booktitle} {AIAA Guidance, Navigation, and Control Conference}}}\ (\bibinfo {year} {2017})\ p.\ \bibinfo {pages} {1734}\BibitemShut {NoStop}%
\bibitem [{\citenamefont {Wu}\ \emph {et~al.}(2009)\citenamefont {Wu}, \citenamefont {Zhang},\ and\ \citenamefont {Sun}}]{wu2009hovering}%
  \BibitemOpen
  \bibfield  {author} {\bibinfo {author} {\bibfnamefont {J.~H.}\ \bibnamefont {Wu}}, \bibinfo {author} {\bibfnamefont {Y.~L.}\ \bibnamefont {Zhang}},\ and\ \bibinfo {author} {\bibfnamefont {M.}~\bibnamefont {Sun}},\ }\bibfield  {title} {\bibinfo {title} {Hovering of model insects: simulation by coupling equations of motion with navier--stokes equations},\ }\href@noop {} {\bibfield  {journal} {\bibinfo  {journal} {Journal of Experimental Biology}\ }\textbf {\bibinfo {volume} {212}},\ \bibinfo {pages} {3313} (\bibinfo {year} {2009})}\BibitemShut {NoStop}%
\bibitem [{\citenamefont {Ellington}(1984)}]{ellington1984aerodynamics}%
  \BibitemOpen
  \bibfield  {author} {\bibinfo {author} {\bibfnamefont {C.~P.}\ \bibnamefont {Ellington}},\ }\bibfield  {title} {\bibinfo {title} {The aerodynamics of hovering insect flight. ii. morphological parameters},\ }\href@noop {} {\bibfield  {journal} {\bibinfo  {journal} {Philosophical Transactions of the Royal Society of London. B, Biological Sciences}\ }\textbf {\bibinfo {volume} {305}},\ \bibinfo {pages} {17} (\bibinfo {year} {1984})}\BibitemShut {NoStop}%
\bibitem [{\citenamefont {Dudley}\ and\ \citenamefont {Ellington}(1990{\natexlab{b}})}]{dudley1990mechanicsbumblebees}%
  \BibitemOpen
  \bibfield  {author} {\bibinfo {author} {\bibfnamefont {R.}~\bibnamefont {Dudley}}\ and\ \bibinfo {author} {\bibfnamefont {C.~P.}\ \bibnamefont {Ellington}},\ }\bibfield  {title} {\bibinfo {title} {Mechanics of forward flight in bumblebees: I. kinematics and morphology},\ }\href@noop {} {\bibfield  {journal} {\bibinfo  {journal} {Journal of Experimental Biology}\ }\textbf {\bibinfo {volume} {148}},\ \bibinfo {pages} {19} (\bibinfo {year} {1990}{\natexlab{b}})}\BibitemShut {NoStop}%
\bibitem [{\citenamefont {Krstic}\ and\ \citenamefont {Wang}(2000)}]{krstic2000stability}%
  \BibitemOpen
  \bibfield  {author} {\bibinfo {author} {\bibfnamefont {M.}~\bibnamefont {Krstic}}\ and\ \bibinfo {author} {\bibfnamefont {H.-H.}\ \bibnamefont {Wang}},\ }\bibfield  {title} {\bibinfo {title} {Stability of extremum seeking feedback for general nonlinear dynamic systems},\ }\href@noop {} {\bibfield  {journal} {\bibinfo  {journal} {Automatica-Kidlington}\ }\textbf {\bibinfo {volume} {36}},\ \bibinfo {pages} {595} (\bibinfo {year} {2000})}\BibitemShut {NoStop}%
\bibitem [{\citenamefont {Bajpai}\ \emph {et~al.}(2024)\citenamefont {Bajpai}, \citenamefont {Elgohary},\ and\ \citenamefont {Eisa}}]{ECC2024}%
  \BibitemOpen
  \bibfield  {author} {\bibinfo {author} {\bibfnamefont {S.}~\bibnamefont {Bajpai}}, \bibinfo {author} {\bibfnamefont {A.~A.}\ \bibnamefont {Elgohary}},\ and\ \bibinfo {author} {\bibfnamefont {S.~A.}\ \bibnamefont {Eisa}},\ }\bibfield  {title} {\bibinfo {title} {Model-free source seeking by a novel single-integrator with attenuating oscillations and better convergence rate: Robotic experiments},\ }in\ \href@noop {} {\emph {\bibinfo {booktitle} {2024 European Control Conference (ECC)}}}\ (\bibinfo {organization} {IEEE},\ \bibinfo {year} {2024})\ pp.\ \bibinfo {pages} {472--479}\BibitemShut {NoStop}%
\bibitem [{\citenamefont {Elgohary}\ \emph {et~al.}(2025)\citenamefont {Elgohary}, \citenamefont {Eisa},\ and\ \citenamefont {Bajpai}}]{unicycle_bio_inspiration}%
  \BibitemOpen
  \bibfield  {author} {\bibinfo {author} {\bibfnamefont {A.~A.}\ \bibnamefont {Elgohary}}, \bibinfo {author} {\bibfnamefont {S.~A.}\ \bibnamefont {Eisa}},\ and\ \bibinfo {author} {\bibfnamefont {S.}~\bibnamefont {Bajpai}},\ }\bibfield  {title} {\bibinfo {title} {Model-free and real-time bioinspired unicycle-based source seeking: Differential wheeled robotic experiments},\ }\href@noop {} {\bibfield  {journal} {\bibinfo  {journal} {arXiv preprint arXiv:2501.02184}\ } (\bibinfo {year} {2025})}\BibitemShut {NoStop}%
\bibitem [{Uni(2025)}]{Unicycle_multiple_positions}%
  \BibitemOpen
  \href {https://www.youtube.com/watch?v=LmipEpx_2aU} {\bibinfo {title} {Source seeking for moving light source: Unicycle-based extremum seeking with gekf}} (\bibinfo {year} {2025})\BibitemShut {NoStop}%
\bibitem [{\citenamefont {Beard}\ \emph {et~al.}(2020)\citenamefont {Beard}, \citenamefont {McLain},\ and\ \citenamefont {Peterson}}]{pid_book}%
  \BibitemOpen
  \bibfield  {author} {\bibinfo {author} {\bibfnamefont {R.~W.}\ \bibnamefont {Beard}}, \bibinfo {author} {\bibfnamefont {T.~W.}\ \bibnamefont {McLain}},\ and\ \bibinfo {author} {\bibfnamefont {C.}~\bibnamefont {Peterson}},\ }\href@noop {} {\emph {\bibinfo {title} {Introduction to Feedback Control: Using Design Studies}}}\ (\bibinfo  {publisher} {Amazon Fulfillment},\ \bibinfo {year} {2020})\BibitemShut {NoStop}%
\bibitem [{\citenamefont {Elgohary}\ and\ \citenamefont {Eisa}(2025{\natexlab{b}})}]{githubdae}%
  \BibitemOpen
  \bibfield  {author} {\bibinfo {author} {\bibfnamefont {A.}~\bibnamefont {Elgohary}}\ and\ \bibinfo {author} {\bibfnamefont {S.}~\bibnamefont {Eisa}},\ }\href@noop {} {\bibinfo {title} {{Hovering-Flight-in-Flapping-with-Extremum-Seeking-Feedback-System}}},\ \bibinfo {howpublished} {\url{https://github.com/MDCL-UC/Hovering-Flapping-with-ES-Feedback-System}} (\bibinfo {year} {2025}{\natexlab{b}}),\ \bibinfo {note} {gitHub}\BibitemShut {NoStop}%
\bibitem [{\citenamefont {Oliveira}\ \emph {et~al.}(2016)\citenamefont {Oliveira}, \citenamefont {Krsti{\'c}},\ and\ \citenamefont {Tsubakino}}]{oliveira2016extremum}%
  \BibitemOpen
  \bibfield  {author} {\bibinfo {author} {\bibfnamefont {T.~R.}\ \bibnamefont {Oliveira}}, \bibinfo {author} {\bibfnamefont {M.}~\bibnamefont {Krsti{\'c}}},\ and\ \bibinfo {author} {\bibfnamefont {D.}~\bibnamefont {Tsubakino}},\ }\bibfield  {title} {\bibinfo {title} {Extremum seeking for static maps with delays},\ }\href@noop {} {\bibfield  {journal} {\bibinfo  {journal} {IEEE Transactions on Automatic Control}\ }\textbf {\bibinfo {volume} {62}},\ \bibinfo {pages} {1911} (\bibinfo {year} {2016})}\BibitemShut {NoStop}%
\bibitem [{\citenamefont {Oliveira}\ and\ \citenamefont {Krstic}(2022)}]{oliveira2022extremum}%
  \BibitemOpen
  \bibfield  {author} {\bibinfo {author} {\bibfnamefont {T.~R.}\ \bibnamefont {Oliveira}}\ and\ \bibinfo {author} {\bibfnamefont {M.}~\bibnamefont {Krstic}},\ }\href@noop {} {\emph {\bibinfo {title} {Extremum seeking through delays and PDEs}}}\ (\bibinfo  {publisher} {SIAM},\ \bibinfo {year} {2022})\BibitemShut {NoStop}%
\bibitem [{\citenamefont {Ru{\v{s}}iti}\ \emph {et~al.}(2018)\citenamefont {Ru{\v{s}}iti}, \citenamefont {Evangelisti}, \citenamefont {Oliveira}, \citenamefont {Gerdts},\ and\ \citenamefont {Krsti{\'c}}}]{ruvsiti2018stochastic}%
  \BibitemOpen
  \bibfield  {author} {\bibinfo {author} {\bibfnamefont {D.}~\bibnamefont {Ru{\v{s}}iti}}, \bibinfo {author} {\bibfnamefont {G.}~\bibnamefont {Evangelisti}}, \bibinfo {author} {\bibfnamefont {T.~R.}\ \bibnamefont {Oliveira}}, \bibinfo {author} {\bibfnamefont {M.}~\bibnamefont {Gerdts}},\ and\ \bibinfo {author} {\bibfnamefont {M.}~\bibnamefont {Krsti{\'c}}},\ }\bibfield  {title} {\bibinfo {title} {Stochastic extremum seeking for dynamic maps with delays},\ }\href@noop {} {\bibfield  {journal} {\bibinfo  {journal} {IEEE control systems letters}\ }\textbf {\bibinfo {volume} {3}},\ \bibinfo {pages} {61} (\bibinfo {year} {2018})}\BibitemShut {NoStop}%
\bibitem [{\citenamefont {Liu}\ and\ \citenamefont {Krstic}(2012)}]{liu2012stochastic}%
  \BibitemOpen
  \bibfield  {author} {\bibinfo {author} {\bibfnamefont {S.-J.}\ \bibnamefont {Liu}}\ and\ \bibinfo {author} {\bibfnamefont {M.}~\bibnamefont {Krstic}},\ }\href@noop {} {\emph {\bibinfo {title} {Stochastic averaging and stochastic extremum seeking}}}\ (\bibinfo  {publisher} {Springer Science \& Business Media},\ \bibinfo {year} {2012})\BibitemShut {NoStop}%
\bibitem [{\citenamefont {Tipsuwan}\ and\ \citenamefont {Chow}(2003)}]{delay}%
  \BibitemOpen
  \bibfield  {author} {\bibinfo {author} {\bibfnamefont {Y.}~\bibnamefont {Tipsuwan}}\ and\ \bibinfo {author} {\bibfnamefont {M.-Y.}\ \bibnamefont {Chow}},\ }\bibfield  {title} {\bibinfo {title} {Control methodologies in networked control systems},\ }\href@noop {} {\bibfield  {journal} {\bibinfo  {journal} {Control engineering practice}\ }\textbf {\bibinfo {volume} {11}},\ \bibinfo {pages} {1099} (\bibinfo {year} {2003})}\BibitemShut {NoStop}%
\bibitem [{\citenamefont {Maggia}\ \emph {et~al.}(2020)\citenamefont {Maggia}, \citenamefont {Eisa},\ and\ \citenamefont {Taha}}]{maggia2020higher}%
  \BibitemOpen
  \bibfield  {author} {\bibinfo {author} {\bibfnamefont {M.}~\bibnamefont {Maggia}}, \bibinfo {author} {\bibfnamefont {S.~A.}\ \bibnamefont {Eisa}},\ and\ \bibinfo {author} {\bibfnamefont {H.~E.}\ \bibnamefont {Taha}},\ }\bibfield  {title} {\bibinfo {title} {On higher-order averaging of time-periodic systems: reconciliation of two averaging techniques},\ }\href@noop {} {\bibfield  {journal} {\bibinfo  {journal} {Nonlinear Dynamics}\ }\textbf {\bibinfo {volume} {99}},\ \bibinfo {pages} {813} (\bibinfo {year} {2020})}\BibitemShut {NoStop}%
\end{thebibliography}%

% \clearpage
\appendix
\section*{Appendix}
The following figures show the state variables' response and the hovering conditions of the natural hovering ES system \eqref{eq:Full_model_ESC} for various insects, including cranefly, bumblebee, dragonfly, hoverfly, and hummingbird. Each figure highlights the behavior of these systems as described in Figure \ref{fig:Hawkmooth_states} and Figure \ref{fig:Hawkmooth_J}.

\begin{figure*}
    \centering
    \includegraphics[width=0.89\linewidth]{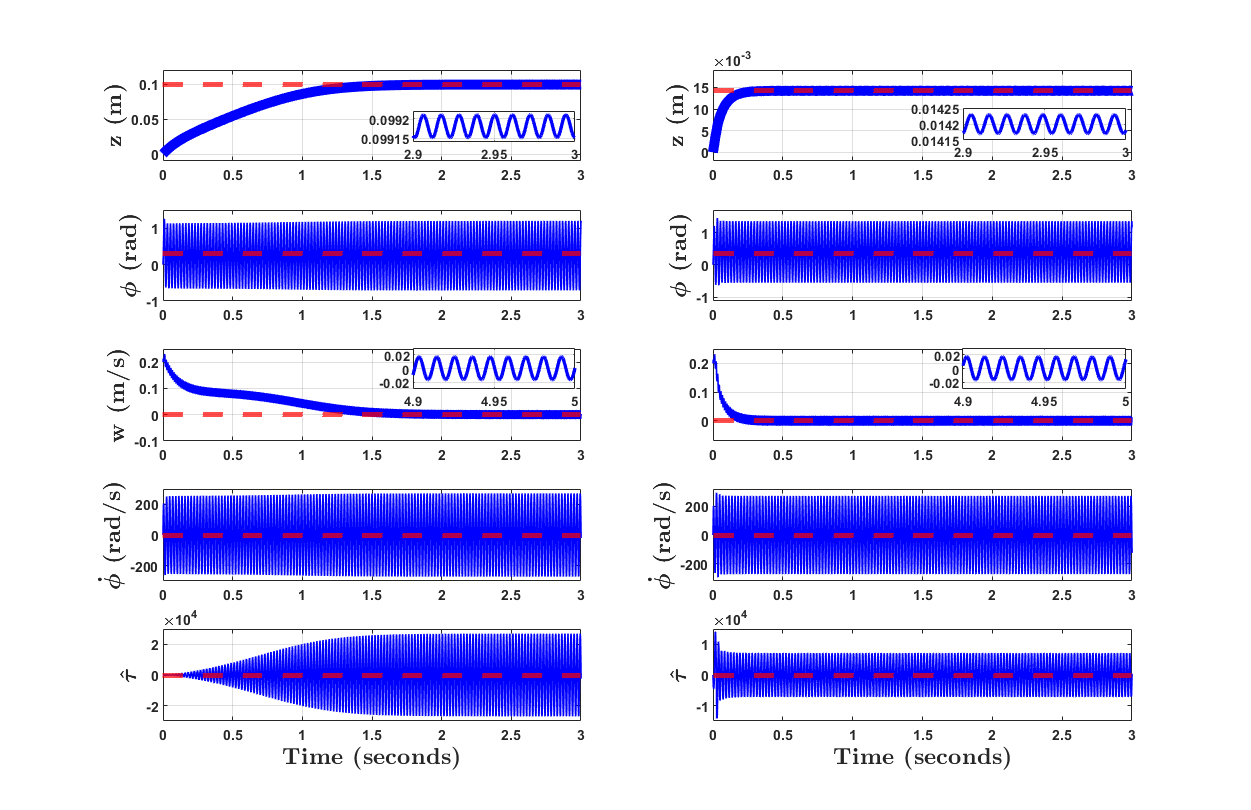}
    \caption{\textcolor{black}{State variables response of the natural hovering ES system \eqref{eq:Full_model_ESC} for \textit{\textbf{cranefly}}. See the caption of Figure \ref{fig:Hawkmooth_states} as this figure follows similar structure.}}
    \label{fig:Cranefly_states}
\end{figure*}

\begin{figure*}
    \centering
\includegraphics[width=0.89\linewidth]{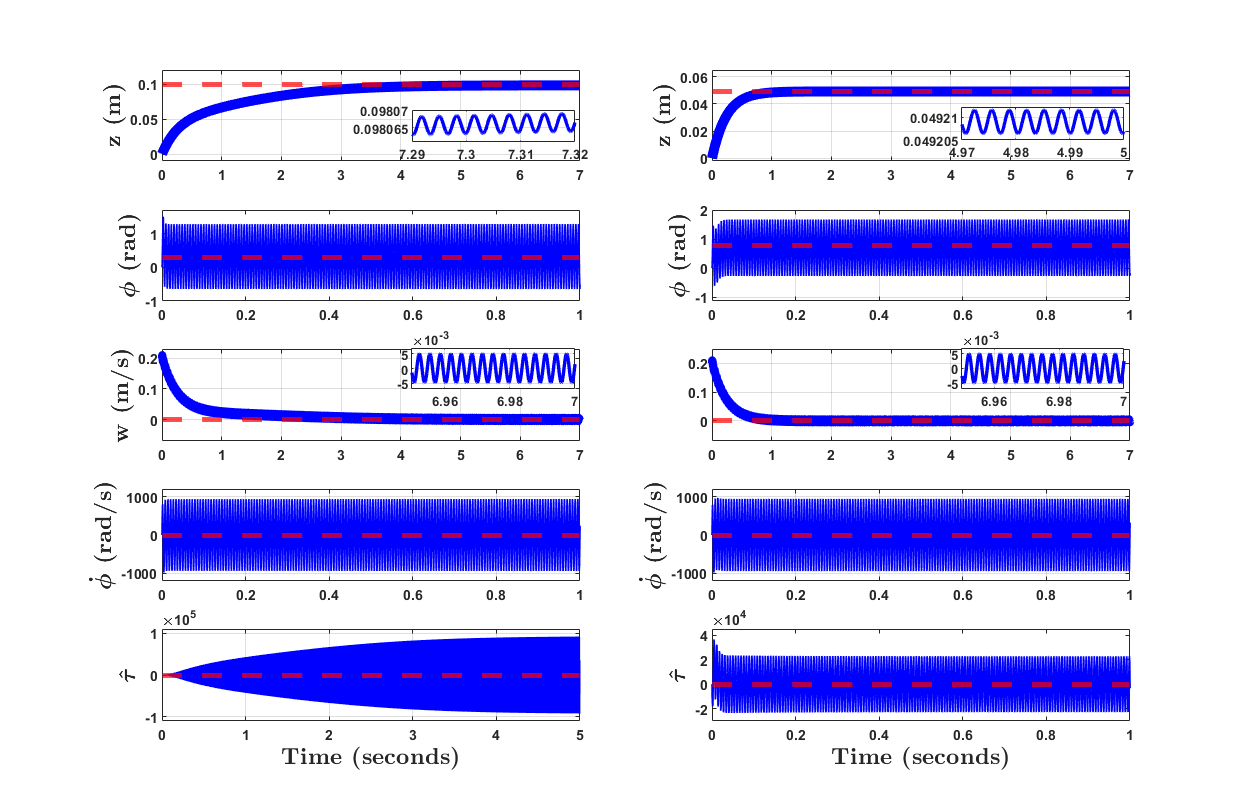}
    \caption{\textcolor{black}{State variables response of the natural hovering ES system \eqref{eq:Full_model_ESC} for \textit{\textbf{bumblebee}}. See the caption of Figure \ref{fig:Hawkmooth_states} as this figure follows similar structure.}}
\label{fig:Bumblebee_states}
\end{figure*}

\begin{figure*}
    \centering
\includegraphics[width=0.89\linewidth]{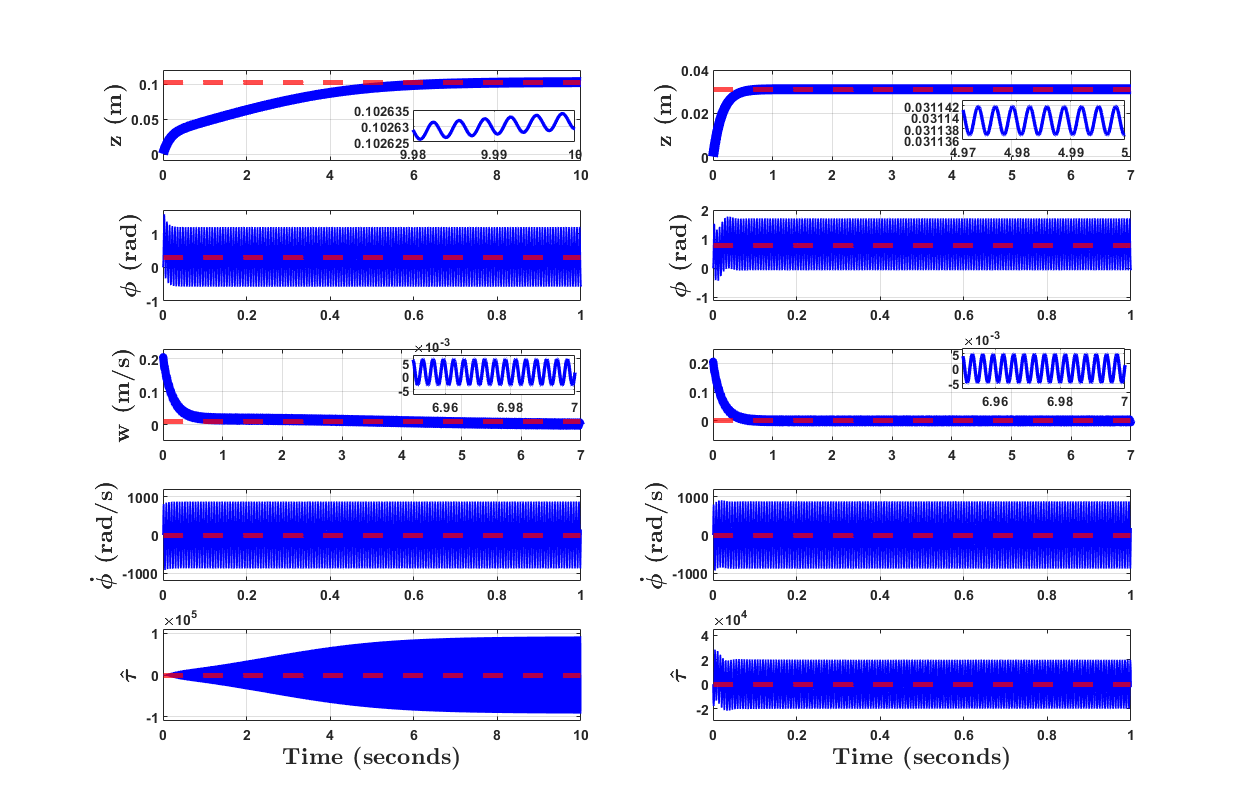}
    \caption{\textcolor{black}{State variables response of the natural hovering ES system \eqref{eq:Full_model_ESC} for \textit{\textbf{dragonfly}}. See the caption of Figure \ref{fig:Hawkmooth_states} as this figure follows similar structure.}}
    \label{fig:Dragonfly_states}
\end{figure*}

\begin{figure*}
    \centering
\includegraphics[width=0.89\linewidth]{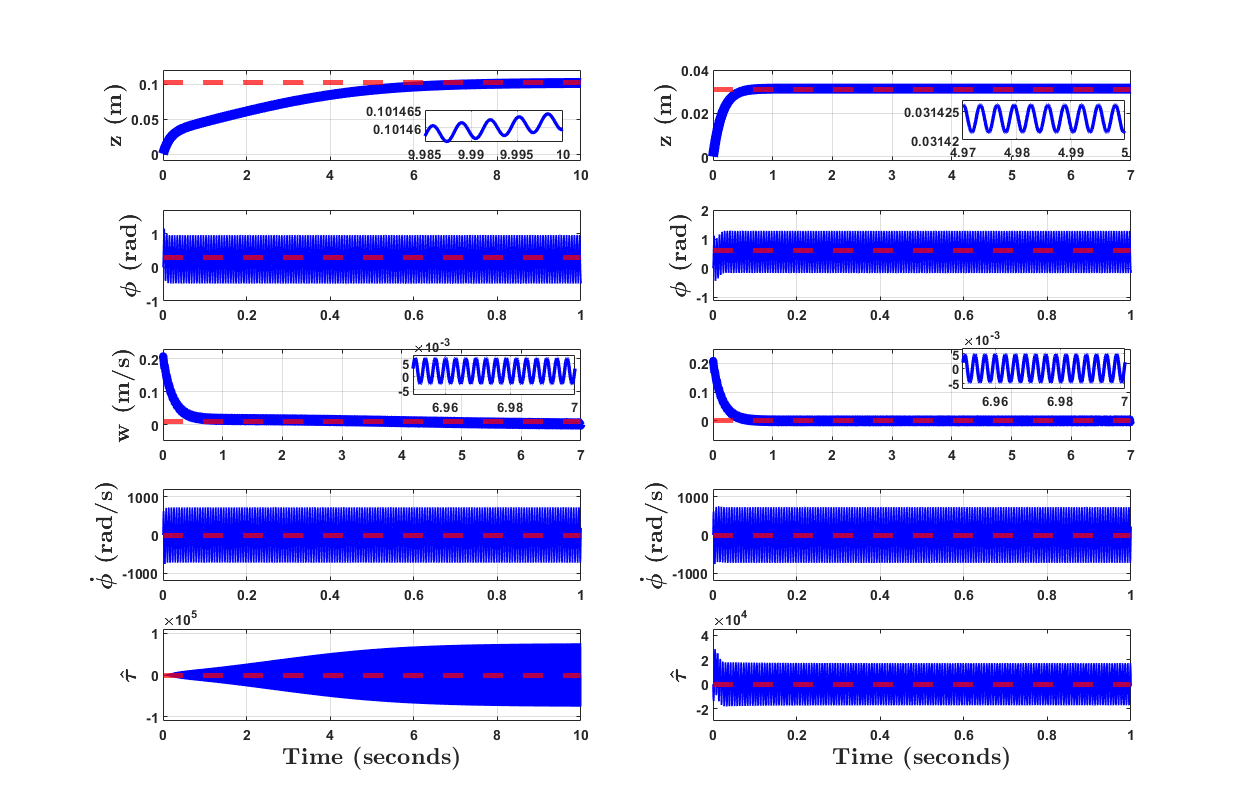}
    \caption{\textcolor{black}{State variables response of the natural hovering ES system \eqref{eq:Full_model_ESC} for \textit{\textbf{hoverfly}}. See the caption of Figure \ref{fig:Hawkmooth_states} as this figure follows similar structure.}}
\label{fig:Hoverfly_states}
\end{figure*}

\begin{figure*}
    \centering
\includegraphics[width=1\linewidth]{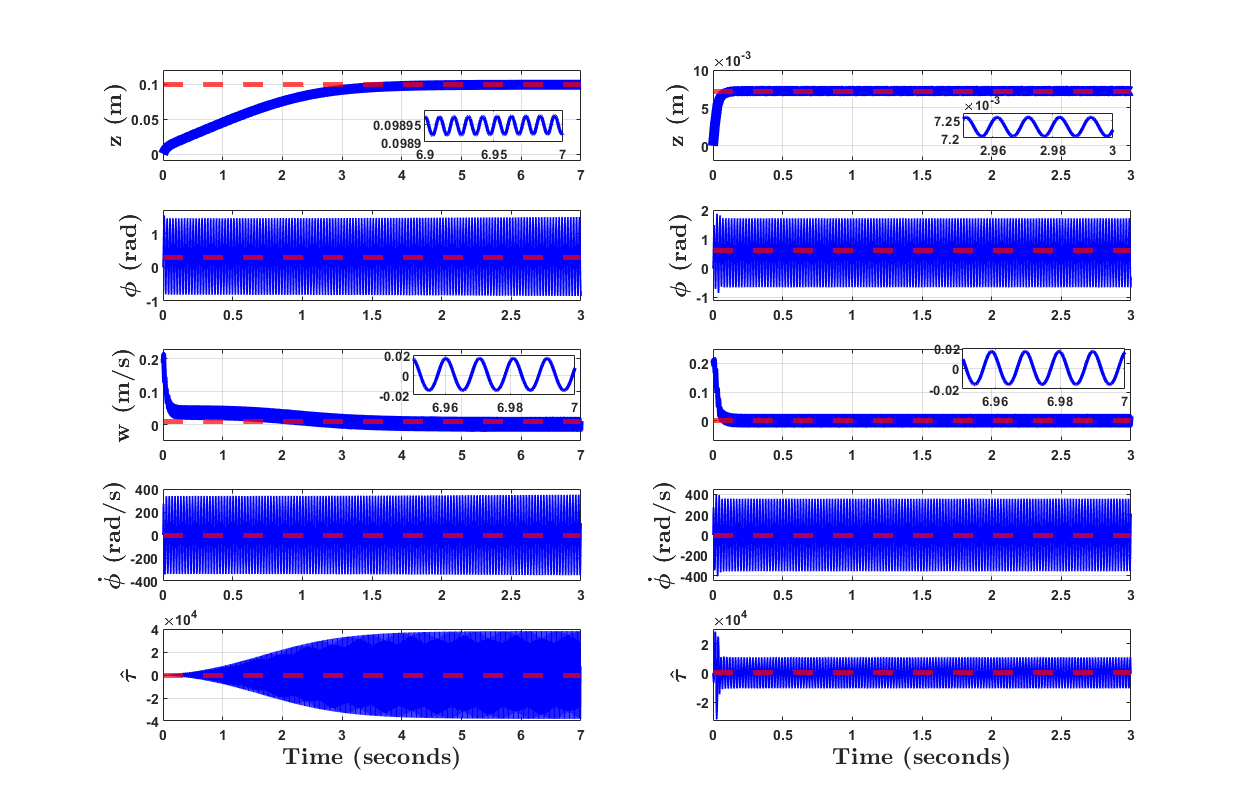}
    \caption{\textcolor{black}{State variables response of the natural hovering ES system \eqref{eq:Full_model_ESC} for \textit{\textbf{a hummingbird}}. See the caption of Figure \ref{fig:Hawkmooth_states} as this figure follows similar structure.}}
\label{fig:Hummingbird_states}
\end{figure*}

\begin{figure*}
    \centering
    \includegraphics[width=1\linewidth]{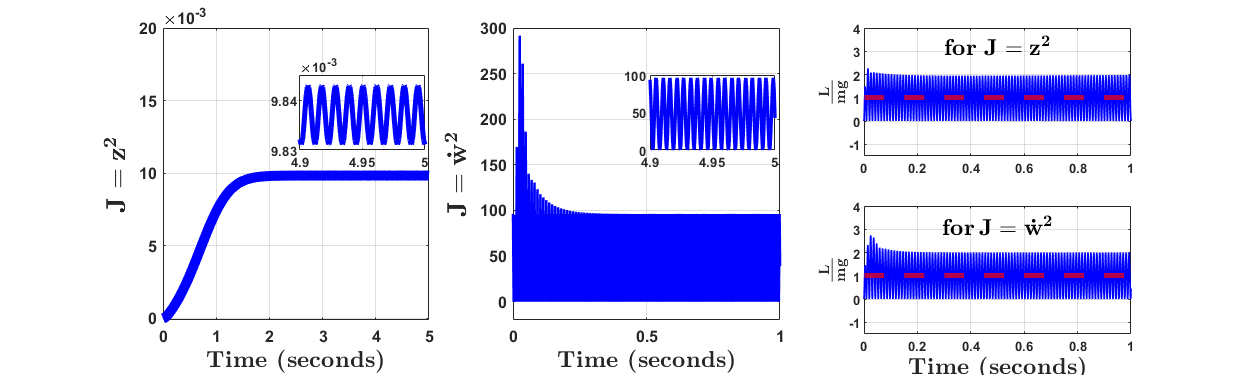}
    \caption{\textcolor{black}{For the \textit{\textbf{cranefly}} case, the hovering condition is observed similar to Figure \ref{fig:Hawkmooth_J}.}}
    \label{fig:Cranefly_J}
\end{figure*}

\begin{figure*}
    \centering
\includegraphics[width=0.9\linewidth]{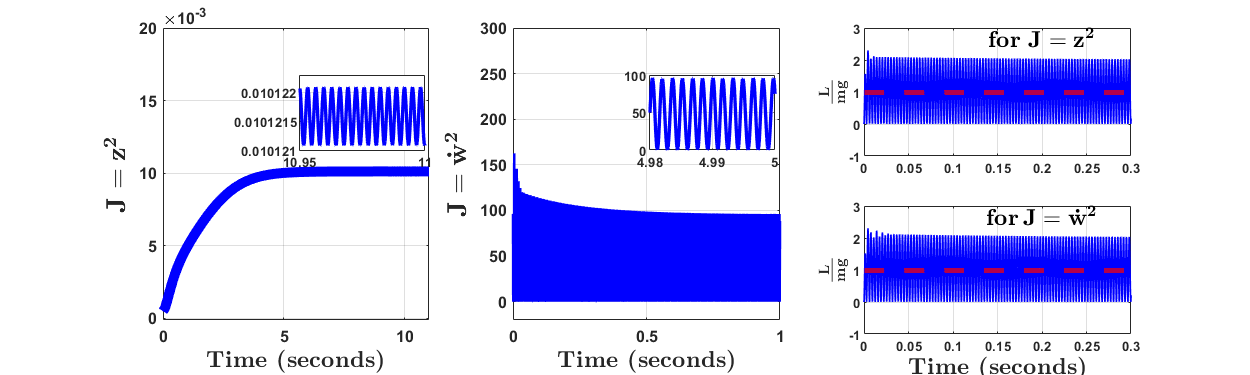}
    \caption{\textcolor{black}{For the \textit{\textbf{bumblebee}} case, the hovering condition is observed similar to Figure \ref{fig:Hawkmooth_J}.}}
    \label{fig:Bumblebee_J}
\end{figure*}

\begin{figure*}
    \centering
    \includegraphics[width=0.9\linewidth]{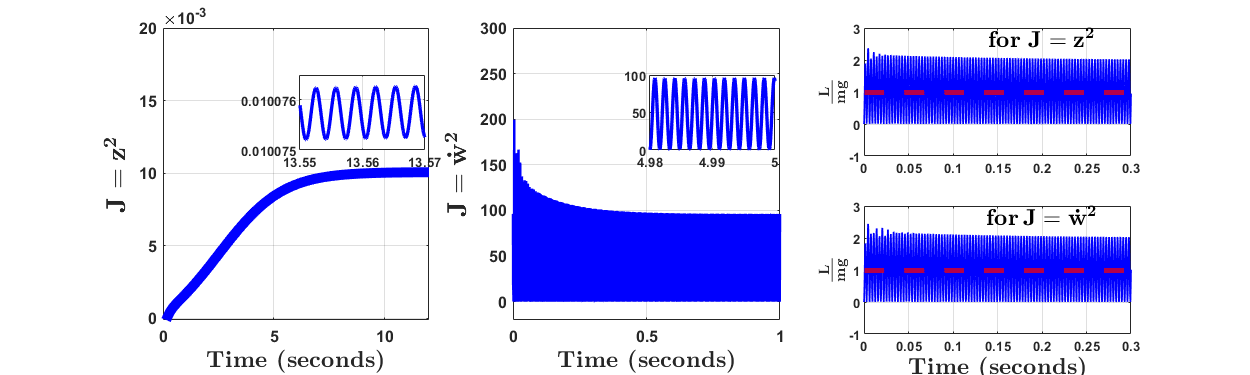}
    \caption{\textcolor{black}{For the \textit{\textbf{dragonfly}} case, the hovering condition is observed similar to Figure \ref{fig:Hawkmooth_J}.}}
    \label{fig:Dragonfly_J}
\end{figure*}

\begin{figure*}
    \centering
\includegraphics[width=0.9\linewidth]{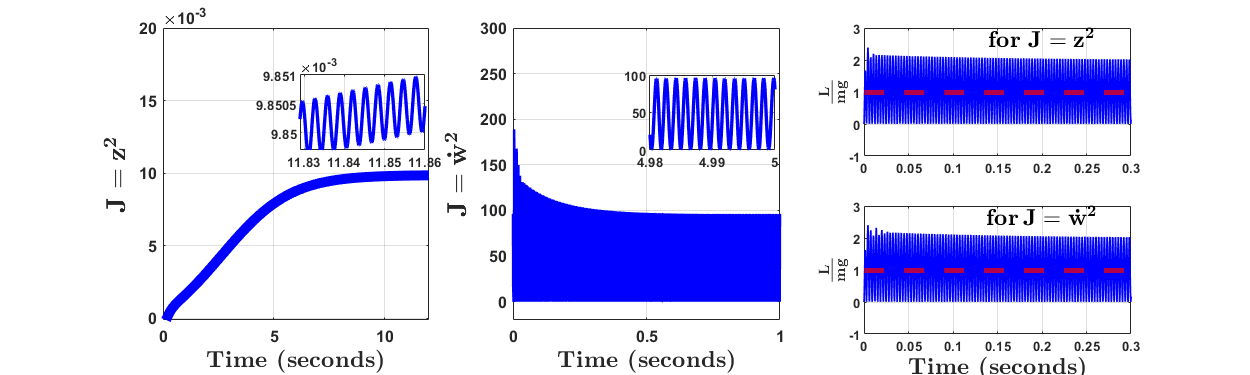}
    \caption{\textcolor{black}{For the \textit{\textbf{hoverfly}} case, the hovering condition is observed similar to Figure \ref{fig:Hawkmooth_J}.}}
    \label{fig:Hoverfly_J}
\end{figure*}

\begin{figure*}
    \centering \includegraphics[width=0.9\linewidth]{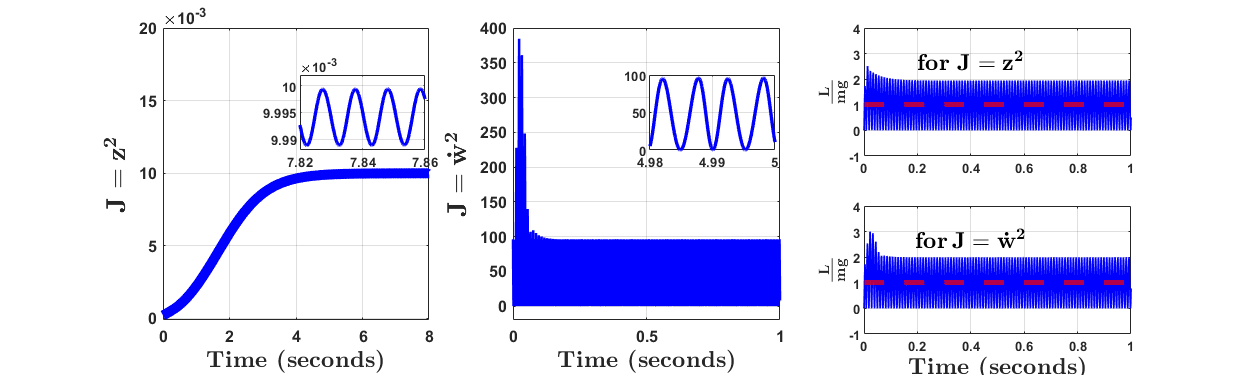}
    \caption{\textcolor{black}{For the \textit{\textbf{hummingbird}} case, the hovering condition is observed similar to Figure \ref{fig:Hawkmooth_J}.}}
    \label{fig:Hummingbird_J}
\end{figure*}

\end{document}